# SCHEDULING A MULTI CLASS QUEUE WITH MANY EXPONENTIAL SERVERS: ASYMPTOTIC OPTIMALITY IN HEAVY TRAFFIC

BY RAMI ATAR[1], AVI MANDELBAUM[2] AND MARTIN I. REIMAN

*Technion–Israel Institute of Technology, Technion–Israel Institute of Technology and Bell Labs, Lucent Technologies*

We consider the problem of scheduling a queueing system in which many statistically identical servers cater to several classes of impatient customers. Service times and impatience clocks are exponential while arrival processes are renewal. Our cost is an expected cumulative discounted function, linear or nonlinear, of appropriately normalized performance measures. As a special case, the cost per unit time can be a function of the number of customers waiting to be served in each class, the number actually being served, the abandonment rate, the delay experienced by customers, the number of idling servers, as well as certain combinations thereof. We study the system in an asymptotic heavy-traffic regime where the number of servers $n$ and the offered load $\mathbf{r}$ are simultaneously scaled up and carefully balanced: $n \approx \mathbf{r} + \beta\sqrt{\mathbf{r}}$ for some scalar $\beta$. This yields an operation that enjoys the benefits of both heavy traffic (high server utilization) and light traffic (high service levels.)

We first consider a formal weak limit, through which our queueing scheduling problem gives rise to a diffusion control problem. We show that the latter has an optimal Markov control policy, and that the corresponding Hamilton–Jacobi–Bellman (HJB) equation has a unique classical solution. The Markov control policy and the HJB equation are then used to define scheduling control policies which we prove are asymptotically optimal for our original queueing system. The analysis yields both qualitative and quantitative insights,

Received August 2002; revised August 2003.

[1]Supported by Israel Science Foundation Grant 126/02, US–Israel Binational Science Foundation Grant 1999179 and the Technion fund for the promotion of research.

[2]Supported by Israel Science Foundation Grants 388/99 and 126/02, by the Niderzaksen Fund and by the Technion funds for the promotion of research and sponsored research.

in particular on staffing levels, the roles of non-preemption and work conservation, and the trade-off between service quality and servers' efficiency.

**Contents**

**1. Introduction.** We analyze a queueing system that consists of several customer classes and a large pool of independent statistically identical servers (see Figure 1). Customer arrivals for each class follow a renewal process. Each server can serve customers of all classes, and service durations are exponentially distributed with class-dependent means. In addition, some customers abandon the system while waiting to be served, and abandonments arise according to exponential clocks with class-dependent rates. This work addresses the stochastic control problem of system scheduling:

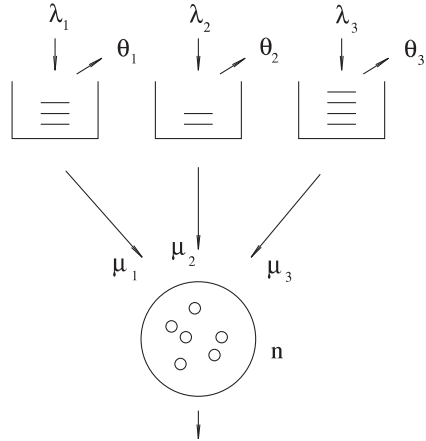

Fig. 1. *A many-server multiclass queueing system.*

how to optimally match customers and servers. The cost criterion we consider is an expected cumulative discounted function of the (appropriately normalized) number of customers waiting to be served and the number actually being served, for each class. Special cases for the cost per unit time are the number of customers in the system (or increasing functions of it), the number of abandonments per unit time, the delay experienced by the customers, the number of idling servers and certain combinations of these costs. Since our scheduling problem is too complex for direct analysis, we resort to heavy-traffic asymptotics. The goal is to identify the asymptotics with a diffusion control problem, then rigorously justify this identification and finally gain insight from it.

1.1. *Motivation: the QED regime.* The asymptotic heavy-traffic regime that we consider is the one analyzed by Jagerman [24], Halfin and Whitt [17] and Fleming, Stolyar and Simon [12]. Here, the number of servers and the arrival rates are large and carefully balanced so that the traffic intensity is moderately close to unity. Economies of scale then enable an operation that is both efficiency-driven (high servers' utilization) and quality-driven (high service levels), hence the terminology QED: both Quality- and Efficiency-driven.

An important motivating application for our model is the modern telephone call center, where a large heterogeneous customer population seeks service from many flexible servers. In this context, the QED regime was identified in practice first by Sze [36], and more recently and systematically in Garnett, Mandelbaum and Reiman [15]. The QED regime captures the operational environment of well-run moderate-to-large call centers, where servers' utilization is high yet a significant fraction of the customers is served immediately upon calling. The last two statements are in fact equivalent for single-class many-server systems [15, 17]. They are further equivalent to "square-root safety staffing," which also applies to the model under study here: if $\mathbf{r}$ denotes the offered load and $n$ the number of servers, then $n \approx \mathbf{r} + \beta\sqrt{\mathbf{r}}$ for some constant $\beta$. (See [13] for more elaboration, motivation and references.)

For a single-class queue $(GI/M/n)$ in the QED regime, one subtracts from the number of customers in the system the number of servers and then divides by the square root of the latter. The resulting stochastic process, when positive, models the (scaled) queue-length, and when negative models the (scaled) number of idle servers. Halfin and Whitt [17] proved that this process converges in distribution, as the number of servers $(n)$ grows without bound, to a diffusion process with a fixed diffusion coefficient and a piecewise linear state-dependent drift, under appropriate assumptions on system parameters. The result was extended in [15] to accommodate abandonment from the queue (but arrivals were assumed Poisson). Further extensions were carried out by Puhalskii and Reiman [34] to cover a multiclass queue, phase-type service time distributions and priority scheduling policies, giving rise in the limit to a multidimensional diffusion process.

1.2. *Diffusion control problems and queueing systems.* There has been a considerable amount of research on diffusion control problems in the context of queueing systems, specifically on asymptotic optimality when approaching a diffusive limit. We refer the reader to [38] for a summary and further references. Most of this research, however, has been within the "conventional" heavy-traffic regime which, in the terminology introduced above, corresponds to an efficiency-driven regime of operation: servers' utilization approaches 100%, with essentially all customers being delayed in queue for service. To wit, our model in "conventional" heavy traffic was analyzed by Van Mieghem [37], who considered a single server (or equivalently, a *fixed* number of servers) with traffic intensity converging to unity. (One could, alternatively, increase the number of servers to infinity, which entails an acceleration of the convergence to unity; see the last section of [30].)

Following Harrison [18], there has been a stream of research that produced schemes for determining "good" scheduling policies for queueing systems, in an asymptotic sense. These have been based on exact analytic solutions to corresponding diffusion control problems, formally obtained as "conventional" heavy traffic limits. For rigorous proofs of asymptotic optimality, see [[3, 26, 27, 29, 30, 31, 32] and [37]].

Recently, Armony and Maglaras [1], Harrison and Zeevi [21] and the present authors [2] have considered stochastic control problems in the QED regime. The first [1] models and analyzes rational customers in equilibrium, and the last [2] served as a pilot for the present paper. The analysis in [21] is that of the diffusion control problem associated with our queueing system with linear costs. Specifically, Harrison and Zeevi show in [21] that this control problem has an optimal Markov control policy (cf. [11]) which is characterized in terms of its underlying HJB equation. Then, they use the diffusion control problem to propose a scheduling control policy for the original queueing system, conjecturing that it is asymptotically optimal in the QED regime. In the current paper we use that same approach, with yet a significant broadening of modeling scope: we identify a sequence of HJB-based scheduling policies (for a general and natural cost structure) and we prove their asymptotic optimality (within a broad family of nonanticipating preemptive or nonpreemptive policies).

1.3. *Main results and scope.* Our main results are as follows. First, we formally take a heavy-traffic limit in the QED regime (Section 2.3). Then we show that the diffusion control problem associated with this limit has an optimal Markov control policy, and that its HJB equation has a unique classical solution (see Theorem 3). This extends the results of [21] to cover a large class of cost functions. As is often the case in stochastic control of diffusions, proving existence of optimal Markov control policies is coupled with establishing the existence and uniqueness of solutions for the underlying HJB equation. In the case of bounded cost, existence and uniqueness for this equation follow from the theory of optimal control of diffusions [6, 11] and of nonlinear elliptic PDEs [23]. Since our cost is not assumed to be bounded, finer information on the model needs to be exploited, and in particular moment estimates on the controlled processes are required [Proposition 4(ii)].

Having studied the diffusion control problem and the HJB equation, we use them to propose a scheme for determining scheduling control policies of two types: preemptive and nonpreemptive (see Section 2.6). After defining a notion of scheduling control policies that do not anticipate the future, we prove that among them, our proposed policies are asymptotically optimal in the QED heavy-traffic limit (Theorems 2 and 4). (More precisely, asymptotic optimality is proved among work conserving policies; more on that in the sequel.) The asymptotic optimality is in the sense that, under the proposed policies, the cost converges to the optimal cost of the diffusion control problem, and that the latter is a lower bound for the limit inferior of costs under any other sequence of policies.

Our approach for deriving the diffusion control problem follows Bell and Williams [3] in that the system of equations and the cost are represented in terms of the system's primitives. The controlled diffusion then arises as a formal weak limit. In obtaining the asymptotic results, this direct relation between the queueing system control problem and the diffusion control problem is convenient.

The policies that we establish as asymptotically optimal are feedback controls. By this we mean that the action at each time depends only on the "state" of the system, namely on the number of customers waiting to be served and the number of customers being served, for each class. The family of policies among which they are proved asymptotically optimal contains all policies that observe all system information up to decision time. In fact, the family we consider is slightly broader in that the policies are allowed to exploit some information on the future, namely the time of the next arrival for each class. We comment below that this is a natural class to consider in the presence of renewal arrivals (cf. Section 2.2).

Under a preemptive scheduling control, service to customers can be interrupted at any time and resumed at a later time. Consequently, the class-fractions of the customers waiting to be served provide natural candidates for control. The diffusion control problem is formulated with such a preemptive model in mind, and the control process corresponds to these fractions (as suggested in [21]). When restricting to scheduling control policies that are nonpreemptive, one must constrain the processes that count the number of customers routed to the server pool to be nondecreasing. The diffusion control problem that arises from such a model resides in a higher dimension. However, here we demonstrate that the nonpreemptive scheduling control problem is asymptotically governed by the simpler diffusion control problem and its HJB equation; to this end, the preemptive HJB equation is used to construct a nonpreemptive scheduling control policy that is asymptotically optimal (in fact, within the class of preemptive policies).

Work-conserving policies are typically not optimal among nonpreemptive scheduling control policies. This can be seen in a simple example, where there are two customer classes, and the cost takes the form $E \int_0^\infty e^{-\gamma t} \sum_{i=1,2} c_i \Phi_i(t) \, dt : \Phi_i(t)$ is the number of class-$i$ customers waiting to be served at time $t$. Consider the event that when the first class-1 customer arrives, there is exactly one free server, and no class-2 waiting customers. If the customer is routed to the free server, then there is a positive probability that the class-2 customer that arrives next will be delayed by at least one unit of time. If the ratio $c_2/c_1$ is large enough, it is clear that the cost paid for delaying this individual class-2 customer can be larger than the cost of delaying all class-1 customers that ever arrive (due to the discount in the cost). As a result, a good policy will leave a free server to idle until a class-2 customer arrives, or until additional servers become idle.

On the other hand, when allowing preemptive policies, for most costs of interest it is intuitively clear that work conservation is optimal. We refer to such costs as *work encouraging* (see Section 5). While there is no attempt here at a rigorous analysis of work encouragement (this seems to require a different modeling framework), our results do reduce the problem of asymptotic optimality (under preemption or nonpreemption) to verifying that work conservation is optimal *among preemptive policies* (Corollary 1). For example, when optimality of work-conserving preemptive policies holds for the prelimit problems, our results, which establish asymptotic optimality of a nonpreemptive policy that is work conserving, imply that the phenomenon described in the previous paragraph is negligible on the diffusive scale.

We comment that, to prove asymptotic optimality, it is not necessary to establish weak convergence of the controlled processes to a controlled diffusion, but only convergence of the costs. However, under appropriate regularity conditions of the coefficients (such as Lipschitz continuity of the function used to define the optimal Markov control policy; see Theorem 3), convergence of the controlled processes follows from our analysis.

Diffusion control problems that arise in "conventional" heavy traffic often have a particularly simple solution, in the form of a static priority policy. Moreover, these policies typically exhibit pathwise minimality of the associated workload processes. Such a simplification is a consequence of a *state-space collapse* [19, 35] namely that these multidimensional diffusion control problems reduce to one-dimensional problems: in conventional heavy traffic, the many servers work in concert as though they constitute a single "super-server." While such collapse prevails in the special case studied in [2], simulations and intuition indicate that, in general for the QED regime, an analogous phenomenon is unlikely to occur. Significantly, though, our analysis does yield some state-space collapse: it is manifested through the asymptotic optimality of nonpreempting work-conserving feedback controls, within the far broader class that allows nonpreemption, idleness in the presence of waiting customers and the use of all past information.

1.4. *Organization and notation.* In Section 2 we describe the model, introduce a notion of scheduling control policies that do not anticipate the future and specify the heavy-traffic assumptions and scaling. We state our first main result regarding the diffusion control problem (Theorem 1). We then use the diffusion control problem to construct two sequences of scheduling control policies (preemptive and nonpreemptive) for the queueing system, and state our second main result on asymptotic optimality of these sequences of policies (Theorem 2). Section 3 treats the diffusion control problem, proving existence and uniqueness for the underlying HJB equation, and existence of optimal Markov control policies. The asymptotic optimality results are proved in Section 4. In Section 5 we discuss the implications of our main result to sequences of policies that are not necessarily work conserving, and chart possible directions for further research. Finally, some auxiliary results are proved in the Appendix.

For $x \in \mathbb{R}^k$ we let $\|x\| = \sum_i |x_i|$. Associated with the parameters $k$ and $n$ of the queueing system are the sets $K = \{1, \ldots, k\}$ and $N = \{1, \ldots, n\}$. We write $\mathbb{N} = \{1, 2, \ldots\}$, $\mathbb{Z}_+^k = \{0, 1, 2, \ldots\}^k$, $\mathbb{R}_+^k = [0, \infty)^k$ and $\mathbb{S}^k = \{x \in \mathbb{R}_+^k : \sum_{i=1}^k x_i = 1\}$.

We denote by $B(m,r)$ an open Euclidean ball of radius $r$ about $m$. $\mathcal{B}(A)$ denotes Borel $\sigma$-field of subsets of $A$. $C^{m,\varepsilon}(D)$ [resp. $C^m(D)$] denotes the class of functions on $D \subset \mathbb{R}^k$ for which all derivatives up to order $m$ are Hölder continuous uniformly on compact subsets of $D$ [continuous on $D$]. $C_{\text{pol}}(\mathbb{R}^k)$ denotes the class of continuous functions $f$ on $\mathbb{R}^k$, satisfying a polynomial growth condition: there are constants $c$ and $r$ such that $|f(x)| \leq c(1 + \|x\|^r)$, $x \in \mathbb{R}^k$. We let $C_{\text{pol}}^{m,\varepsilon} = C_{\text{pol}} \cap C^{m,\varepsilon}$. For $E$ a metric space, we denote by $\mathbb{D}(E)$ the space of all cadlag functions (i.e., right continuous and having left limits) from $\mathbb{R}_+$ to $E$. We endow $\mathbb{D}(E)$ with the usual Skorohod topology. All processes we consider are assumed to have sample paths in $\mathbb{D}(E)$ (for appropriate $E$, mostly $E = \mathbb{R}^k$). If $X^n$, $n \in \mathbb{N}$ and $X$ are processes with sample paths in $\mathbb{D}(E)$, we write $X^n \Rightarrow X$ to denote weak convergence of the measures induced by $X^n$ [on $\mathbb{D}(E)$] to the measure induced by $X$. For any cadlag path $X$, let $X_{t-} = \lim_{s \uparrow t} X_s$ for $t > 0$, $X_{0-} = X_0$, and $\Delta X_t = X_t - X_{t-}$. If $X$ is a process (or a function on $\mathbb{R}_+$), $\|X\|_t^* = \sup_{0 \leq s \leq t} \|X(s)\|$, and if $X$ takes real values, $|X|_t^* = \sup_{0 \leq s \leq t} |X(s)|$. $X(t)$ and $X_t$ are used interchangeably. Vectors in $\mathbb{R}^k$ are considered as column vectors. We write $\mathbb{1} = (1,\ldots,1)' \in \mathbb{R}^k$. For vectors $u, v \in \mathbb{R}^k$, let $u \cdot v$ denote their scalar product. Finally, $c$ denotes a positive constant whose value is not important, and may change from line to line.

**2. The controlled system in the QED regime and its diffusion approximation.** We consider a queueing system which consists of $k$ customer classes and $n$ multiskilled servers (see Figure 1). Service to any customer can be provided by any of the servers indifferently. The service time distribution depends on the customer class, but not on the individual server (or customer). We say that a customer is in queue $i$ at time $t$ if the customer is of class $i$, and at time $t$ it is in the system and is not being served (although it possibly received partial service prior to time $t$). Customers enter the system at one of the queues, and leave the system in one of two ways: either when their service is completed, or while they are waiting at their queue and decide to abandon the system without being served.

2.1. *The stochastic model.* Let a complete probability space, $(\Omega, F, P)$ be given, on which all the stochastic processes below are defined. Expectation with respect to $P$ is denoted by $E$. The parameter $n$, denoting the number of servers, which is particularly significant in our analysis, will appear (as a superscript) in the notation of all basic stochastic processes associated with the queueing system.

For $i \in K$, the number of class-$i$ customers in the queue at time $t \geq 0$ is denoted by $\Phi_i^n(t)$, and $\Phi^n(t) = (\Phi_1^n(t),\ldots,\Phi_k^n(t))'$. The number of class-$i$ customers being served at time $t$ is denoted by $\Psi_i^n(t)$ and $\Psi^n(t) = (\Psi_1^n(t),\ldots,\Psi_k^n(t))'$.

Clearly these processes take integer values, and

(1) $$\Phi^n(t), \Psi^n(t) \in \mathbb{R}_+^k, \qquad \sum_i \Psi_i^n(t) \le n, \qquad t \ge 0.$$

The initial conditions of the system are assumed to be deterministic and are denoted by $\Phi^n(0) = \Phi^{0,n} = (\Phi_1^{0,1}, \ldots, \Phi_k^{0,1})'$ and $\Psi^n(0) = \Psi^{0,n} = (\Psi_1^{0,n}, \ldots, \Psi_k^{0,n})'$.

Let $A_i^n$, $i \in K$, be independent renewal processes defined as follows. For $i \in K$, let there be a sequence $\{\check{U}_i(j), j \in \mathbb{N}\}$ of strictly positive i.i.d. random variables with mean $E\check{U}_i(1) = 1$ and squared coefficient of variation $\operatorname{Var}(\check{U}_i(1))/(E\check{U}_i(1))^2 = C_{U,i}^2 \in [0, \infty)$. Let

(2) $$U_i^n(j) = \frac{1}{\lambda_i^n} \check{U}_i(j), \qquad i \in K, \ j \in \mathbb{N},$$

where $\lambda_i^n > 0$. With the convention $\sum_1^0 = 0$, define

(3) $$A_i^n(t) = \sup\left\{ m \ge 0 : \sum_{j=1}^m U_i^n(j) \le t \right\}, \qquad i \in K, \ t \ge 0.$$

The value $A_i^n(t)$ denotes the number of arrivals of class-$i$ customers up to time $t$. Note that the first class-$i$ customer arrives at $U_i^n(1)$, and the time between the $(m-1)$st and $m$th arrival of class-$i$ customers is $U_i^n(m)$, $m = 2, 3, \ldots$.

The service time of a class-$i$ customer is assumed to be exponentially distributed with parameter $\mu_i^n$, regardless of the service provider. This is captured in the following description. For $i \in K$, let $S_i^n$ be a Poisson process of rate $\mu_i^n \in (0, \infty)$, and assume that the processes $S_i^n$ are independent of each other and of the processes $A_i^n$, $i \in K$. Let $T_i^n(t)$ denote the time up to $t$ that a server has devoted to class-$i$ customers, summed over all servers. Clearly,

$$T_i^n(t) = \int_0^t \Psi_i^n(s)\, ds, \qquad i \in K, \ t \ge 0.$$

Then $S_i^n(T_i^n(t)) = S_i^n(\int_0^t \Psi_i^n(s)\, ds)$ denotes the number of service completions of class-$i$ jobs, by all servers, up to time $t$. Our assumptions on $T^n$ will ensure that, for each $t$, $T^n(t)$ is independent of any increment of the form $S^n(T^n(t) + s) - S^n(T^n(t))$, $s \ge 0$ (cf. Definition 2).

For $i \in K$, individuals abandon queue $i$ at rate $\theta_i^n \in [0, \infty)$. Let $R_i^n$ be Poisson processes of rate $\theta_i^n$, independent of each other and of the processes $A_j^n, S_j^n$, $j \in K$. Note that the time up to $t$ that a class-$i$ customer spends in the queue, summed over all customers, is equal to $\int_0^t \Phi_i^n(s)\, ds$. Then $R_i^n(\int_0^t \Phi_i^n(s)\, ds)$ denotes the number of abandonments from queue $i$ up to time $t$. Under an appropriate assumption on $\int_0^\cdot \Phi^n(s)\,ds$, similar to that on $T^n$ (cf. Definition 2), this describes abandonment of class-$i$ customers according to independent rate-$\theta_i^n$ Poisson clocks, each run as long as the customer is in the queue.

We would like to have equations that hold for both nonpreemptive and preemptive resume policies. Consider the processes $B_i^n(t)$, $i \in K$, described as follows. $B_i^n(0) = 0$; $B_i^n$ increases by 1 each time a class-$i$ job is assigned to a server (to start or resume service), and decreases by 1 each time such a job is moved back to the queue (in a preemptive-resume policy). Note that in a nonpreemptive policy, $B_i^n(t)$ is the number of type-$i$ customers that have been routed to the server pool at any time up to $t$. In fact, we do not assume that these processes only jump by $\pm 1$; their increments can take arbitrary values in $\mathbb{Z}$. Following are the system equations:

$$\Phi_i^n(t) = \Phi_i^{0,n} + A_i^n(t) - B_i^n(t) - R_i^n\left(\int_0^t \Phi_i^n(s)\,ds\right), \qquad i \in K,\ t \geq 0,$$
(4)
$$\Psi_i^n(t) = \Psi_i^{0,n} + B_i^n(t) - S_i^n\left(\int_0^t \Psi_i^n(s)\,ds\right), \qquad i \in K,\ t \geq 0.$$

These equations hold regardless of assumptions on the policy as to whether it is preemptive or not, and work conserving or not (these terms are, in fact, made precise later in this section). Note that the representations above in terms of Poisson processes $S_i^n$ and $R_i^n$ exploit the exponential assumptions on service times and abandonment.

Assume that there is a full $P$-measure set under which all $A_i^n(t) < \infty$ for $t \geq 0$, $A_i^n$ increases to infinity, $\Delta A_i^n(t) \in \{0,1\}$ for all $t$, and where similar statements hold for $S_i^n$ and $R_i^n$. Then, without loss, we omit from subsequent discussions all realizations (sample paths) of these processes that do not adhere to these conditions.

Let
$$X^n(t) = \Phi^n(t) + \Psi^n(t) \tag{5}$$

and denote $X^{0,n} = \Phi^{0,n} + \Psi^{0,n}$. Then $X_i^n(t)$ is equal to the number of class-$i$ customers in the system at time $t$. The constraints (1) can be written in terms of $X^n$ and $\Psi^n$ as

$$(6)\quad X^n(t) - \Psi^n(t) \in \mathbb{R}_+^k, \qquad \Psi^n(t) \in \mathbb{R}_+^k, \qquad \sum_i \Psi_i^n(t) \leq n, \qquad t \geq 0,$$

while the system equations (4) imply that

$$X_i^n(t) = X_i^{0,n} + A_i^n(t) - R_i^n\left(\int_0^t (X_i^n(s) - \Psi_i^n(s))\,ds\right) - S_i^n\left(\int_0^t \Psi_i^n(s)\,ds\right),$$
(7)
$$i \in K,\ t \geq 0.$$

2.2. *Scheduling control policies.* We define two types of control problems, one where scheduling is preemptive and one where it is nonpreemptive. Equation (7) serves as the description of the system dynamics. The scheduling control policy (SCP) will be identified with the process $\Psi^n$, and it will be assumed that it is such that the constraints (6) are satisfied. Apart from a nonanticipating assumption on $\Psi^n$ (Definition 2), there will be no further restrictions for preemptive scheduling control problems. For nonpreemptive scheduling control problems, a further constraint will be that the process $B^n$ is nondecreasing in each component.

For the following definition, note that, given a process $\Psi^n$, if there exists a process $X^n$ so that (7) holds, then it is unique (as can be argued by induction on the jump times of the processes $A^n$, $R^n$ and $S^n$). Thus (5) uniquely determines $\Phi^n$, and either part of (4) then uniquely determines $B^n$. Also, finiteness of the integrals appearing in (4) and (7) follows from the fact that $\Psi^n_i$ are bounded by $n$, while $X^n_i(t) - \Psi^n_i(t) = \Phi^n_i(t) \leq X^{0,n} + A^n_i(t)$.

DEFINITION 1. (i) We say that $\Psi^n$ is a *preemptive resume scheduling control policy* (P-SCP) if it is a stochastic process with cadlag paths, taking values in $\mathbb{R}^k$, for which there exists a process $X^n$ (referred to as a *controlled process*) satisfying the system equations (7), and such that the constraints (6) are met. Given a P-SCP $\Psi^n$ and a controlled process $X^n$, denote by $\Phi^n$ and $B^n$ the processes uniquely determined by (4) and (5).

(ii) We say that $\Psi^n$ is a *nonpreemptive scheduling control policy* (N-SCP) if it is a P-SCP, and in addition, $B^n_i$, $i \in K$, have nondecreasing paths.

We collectively refer to P-SCPs and N-SCPs as *scheduling control policies* (SCPs) (although the class of SCPs is simply the class of P-SCPs).

We need a notion of SCPs that do not anticipate the future. To this end, denote

$$(8) \qquad T^n_i(t) = \int_0^t \Psi^n_i(s)\,ds, \qquad \mathring{T}^n_i(t) = \int_0^t \Phi^n_i(s)\,ds,$$

and for $i \in K$, let

$$\tau^n_i(t) = \inf\{u \geq t : A^n_i(u) - A^n_i(u-) > 0\}$$

stand for the time of the first arrival to queue $i$ no earlier than $t$. Set

$$(9) \quad \mathcal{F}^n_t = \sigma\{A^n_i(s), S^n_i(T^n_i(s)), R^n_i(\mathring{T}^n_i(s)), \Phi^n_i(s), \Psi^n_i(s), X^n_i(s) : i \in K, s \leq t\}$$

and

$$\begin{aligned}(10)\quad \mathcal{G}^n_t = \sigma\{&A^n_i(\tau^n_i(t) + u) - A^n_i(\tau^n_i(t)), S^n_i(T^n_i(t) + u) - S^n_i(T^n_i(t)),\\ &R^n_i(\mathring{T}^n_i(t) + u) - R^n_i(\mathring{T}^n_i(t)) : i \in K, u \geq 0\}.\end{aligned}$$

While $\mathcal{F}_t^n$ represents the information available at time $t$, $\mathcal{G}_t^n$ constitutes future information. Since for each $i$, $A_i^n$ is a renewal process, its increments of the form that appears in the definition of $\mathcal{G}_t^n$ are independent of $\sigma\{A_i^n(s): s \leq t\}$. However, the time $\tau_i^n$ of the next arrival may be anticipated, to some degree, from the information on the arrivals up to time $t$. Therefore, with $\tau_i^n(t)$ replaced by $t$ in its definition, $\mathcal{G}_t^n$ would not be a good candidate to represent innovative information. Note that an analogous treatment of $S^n$ and $R^n$ is not necessary, since these are Poisson processes which are memoryless. The following definition refers to both types of problems.

DEFINITION 2. We say that a scheduling control policy is *admissible* if:

(i) for each $t$, $\mathcal{F}_t^n$ is independent of $\mathcal{G}_t^n$;
(ii) for each $i$ and $t$, the process $S_i^n(T_i^n(t) + \cdot) - S_i^n(T_i^n(t))$ is equal in law to $S_i^n(\cdot)$, and the process $R_i^n(\mathring{T}_i^n(t) + \cdot) - R_i^n(\mathring{T}_i^n(t))$ is equal in law to $R_i^n(\cdot)$.

Some SCPs considered in this paper will be constructed by setting

$$(11) \qquad \Psi^n(t) = F(X^n(t)), \qquad t \geq 0,$$

for an appropriate choice of $F$. As the following result shows, this leads to admissible SCPs.

PROPOSITION 1. *Fix $n$ and let a function $F:\mathbb{Z}_+^k \to \mathbb{Z}_+^k$ be given such that, for $X \in \mathbb{Z}_+^k$, one has $X - F(X) \in \mathbb{Z}_+^k$ and $\mathbb{1} \cdot F(X) \leq n$. Then the system of equations (7) and (11) has a unique solution, and $\Psi^n$ is an admissible SCP. In particular, if the process $B^n$ determined via (4) has nondecreasing paths, $\Psi^n$ is an admissible N-SCP.*

See the Appendix for a proof.

2.3. *QED scaling.* We consider a sequence of queueing systems as above where now the number of servers $n \in \mathbb{N}$ is used as an index to the sequence. It is implicitly assumed that there is an SCP associated with each queueing system. It is assumed (without loss) that there is one probability space, $(\Omega, F, P)$, on which the processes associated with the $n$th system are defined, for all $n \in \mathbb{N}$. The heavy-traffic assumptions are as follows (cf. [15, 17, 21, 34]).

ASSUMPTION 1. (i) *Parameters.* There are constants $\lambda_i, \mu_i \in (0, \infty)$, $\theta_i \in [0, \infty)$, $\hat{\lambda}_i, \hat{\mu}_i \in \mathbb{R}$, $i \in K$, such that

$$\sum_{i=1}^k \lambda_i/\mu_i = 1$$

and, as $n \to \infty$,

$$n^{-1}\lambda_i^n \to \lambda_i, \qquad \mu_i^n \to \mu_i, \qquad \theta_i^n \to \theta_i$$
$$n^{1/2}(n^{-1}\lambda_i^n - \lambda_i) \to \hat{\lambda}_i, \qquad n^{1/2}(\mu_i^n - \mu_i) \to \hat{\mu}_i.$$

(ii) *Initial conditions.* There are constants $\phi_i \in [0, \infty)$, $\psi_i \in \mathbb{R}$, $i \in K$, such that $\sum_K \psi_i \leq 0$, and, with $\rho_i = \lambda_i/\mu_i$, as $n \to \infty$,

$$\hat{\Phi}_i^{0,n} := n^{-1/2}\Phi_i^{0,n} \to \phi_i, \qquad \hat{\Psi}_i^{0,n} := n^{-1/2}(\Psi_i^{0,n} - \rho_i n) \to \psi_i.$$

REMARK 1. The above scaling is in concert with that in [12, 15, 17, 24]. For a verification, let $\rho^n$ denote the traffic intensity of our $n$th system. Then $\rho^n = \mathbf{r}^n/n$, where its offered load $\mathbf{r}^n$ is given by

$$\mathbf{r}^n = \sum_{i=1}^k \lambda_i^n/\mu_i^n.$$

From Assumption 1 it now follows, via simple algebra, that

$$\sqrt{n}(1 - \rho^n) \to \sum_{i=1}^k (\rho_i\hat{\mu}_i - \hat{\lambda}_i)/\mu_i.$$

Denoting this last limit by $\beta$, we deduce that

$$n \approx \mathbf{r}^n + \beta\sqrt{\mathbf{r}^n}.$$

QED scaling thus leads to square-root safety staffing [7], which characterizes the regimes in [12, 15, 17, 24]. ($\beta > 0$ was required in the original Halfin–Whitt regime of [17], to guarantee stability when there is no abandonment. Our analysis, however, covers all values of $\beta$ since it does not require stability of the queueing system. Indeed, the total discounted costs are always finite in view of our polynomial growth constraints on the cost functions.)

For more details on QED scaling, readers are referred to [15] and [17]. An instructive comparison of the QED regime with conventional heavy traffic, in the context of our problem, is provided by [21].

The rescaled processes are defined as follows:

$$\bar{\Phi}_i^n(t) = n^{-1}\Phi_i^n(t), \qquad \bar{\Psi}_i^n(t) = n^{-1}\Psi_i^n(t),$$
$$\bar{X}_i^n(t) := \bar{\Phi}_i^n(t) + \bar{\Psi}_i^n(t) = n^{-1}X_i^n(t),$$
$$\hat{\Phi}_i^n(t) = n^{1/2}\bar{\Phi}_i^n(t) = n^{-1/2}\Phi_i^n(t),$$
$$\hat{\Psi}_i^n(t) = n^{1/2}(\bar{\Psi}_i^n(t) - \rho_i) = n^{-1/2}(\Psi_i^n(t) - \rho_i n),$$
$$\hat{X}_i^n(t) := \hat{\Phi}_i^n(t) + \hat{\Psi}_i^n(t) = n^{1/2}(\bar{X}_i^n(t) - \rho_i) = n^{-1/2}(X_i^n(t) - \rho_i n).$$

The primitive processes are rescaled as
$$\hat{A}_i^n(t) = n^{-1/2}(A_i^n(t) - \lambda_i^n t), \qquad \hat{S}_i^n(t) = n^{-1/2}(S_i^n(nt) - n\mu_i^n t),$$
$$\hat{R}_i^n(t) = n^{-1/2}(R_i^n(nt) - n\theta_i^n t).$$
Finally,
$$\hat{B}_i^n(t) = n^{-1/2}(B_i^n(t) - n\lambda_i t).$$
With this notation, the system equations (4) can be written as follows:
$$\hat{\Phi}_i^n(t) = \hat{\Phi}_i^{0,n} + \hat{A}_i^n(t) + n^{1/2}(n^{-1}\lambda_i^n - \lambda_i)t$$
$$- \hat{B}_i^n(t) - \hat{R}_i^n\bigg(\int_0^t \bar{\Phi}_i^n(s)\,ds\bigg) - \theta_i^n \int_0^t \hat{\Phi}_i^n(s)\,ds,$$
(12)
$$\hat{\Psi}_i^n(t) = \hat{\Psi}_i^{0,n} + \hat{B}_i^n(t) - \hat{S}_i^n\bigg(\int_0^t \bar{\Psi}_i^n(s)\,ds\bigg)$$
$$- \mu_i^n \int_0^t \hat{\Psi}_i^n(s)\,ds - \rho_i n^{1/2}(\mu_i^n - \mu_i)t.$$

We have from (12)
$$(13) \quad \hat{X}_i^n(t) = \hat{X}_i^{0,n} + r_i \hat{W}_i^n(t) + \ell_i^n t - \mu_i^n \int_0^t \hat{\Psi}_i^n(s)\,ds - \theta_i^n \int_0^t \hat{\Phi}_i^n(s)\,ds,$$
where we denote
$$r_i \hat{W}_i^n(t) = \hat{A}_i^n(t) - \hat{S}_i^n\bigg(\int_0^t \bar{\Psi}_i^n(s)\,ds\bigg) - \hat{R}_i^n\bigg(\int_0^t \bar{\Phi}_i^n(s)\,ds\bigg),$$
(14)
$$r_i = (\lambda_i C_{U,i}^2 + \lambda_i)^{1/2}$$
and
$$\ell_i^n = n^{1/2}(n^{-1}\lambda_i^n - \lambda_i) - \rho_i n^{1/2}(\mu_i^n - \mu_i).$$

We now present a formal derivation of the limiting dynamics, as described by a system of controlled SDEs. The actual relation to the sequence of queueing systems (as a limit) will be justified once our results of Section 4 are established. To this end, we pretend that the convergence
$$\bar{\Phi}_i^n(\cdot) \Rightarrow 0, \qquad \bar{\Psi}_i^n(\cdot) \Rightarrow \rho_i,$$
holds, and write $A, S, \Phi, \Psi, X, B$ for the formal weak limits of $\hat{A}^n, \hat{S}^n, \hat{\Phi}^n, \hat{\Psi}^n$, $\hat{X}^n, \hat{B}^n$ (without worrying at this point about whether weak limits exist). For $i \in K$, the processes $A_i$ and $S_i$ are Brownian motions with zero drift and variances $\lambda_i C_{U,i}^2$ and $\mu_i$, respectively. We thus obtain
$$\Phi_i(t) = \phi_i + A_i(t) + \hat{\lambda}_i t - B_i(t) - \theta_i \int_0^t \Phi_i(s)\,ds,$$
(15)
$$\Psi_i(t) = \psi_i + B_i(t) - \rho_i^{1/2} S_i(t) - \mu_i \int_0^t \Psi_i(s)\,ds - \rho_i \hat{\mu}_i t.$$

The corresponding constraints are as follows:
$$\Phi_i(t) \geq 0, \qquad \sum_i \Psi_i(t) \leq 0.$$

Writing $W = (W_1, \ldots, W_k)'$, $W_i = r_i^{-1}(A_i - \rho_i^{1/2} S_i)$, the process $W$ is a standard $k$-dimensional Brownian motion. The process $X = \Phi + \Psi$ then satisfies

(16) $\quad X_i(t) = x_i + r_i W_i(t) + \ell_i t - \theta_i \int_0^t (X_i(s) - \Psi_i(s))\, ds - \mu_i \int_0^t \Psi_i(s)\, ds,$

as well as the constraints

(17) $\qquad X_i(t) - \Psi_i(t) \geq 0, \qquad \sum_i \Psi_i(t) \leq 0,$

where
$$\ell_i = \hat{\lambda}_i - \rho_i \hat{\mu}_i, \qquad x_i = \phi_i + \psi_i.$$

2.4. *Work conservation and cost.* A policy is work conserving if there can be no idling servers when there are customers in the queue. For the following definition, recall that $\mathbb{1} \cdot \Phi^n$ equals the number of customers in all queues, and that $\mathbb{1} \cdot X^n$ equals the number of customers in the system.

DEFINITION 3. We say that an SCP is *work-conserving* if

(18) $\qquad (\mathbb{1} \cdot X^n(t) - n)^+ = \mathbb{1} \cdot \Phi^n(t), \qquad t \geq 0.$

Note that equivalently

(19) $\qquad (\mathbb{1} \cdot \hat{X}^n(t))^+ = \mathbb{1} \cdot \hat{\Phi}^n(t), \qquad t \geq 0.$

For a given SCP, let $\hat{\Phi}^n$ and $\hat{\Psi}^n$ denote the rescaled processes as before. We consider the problem of infimizing an expected cumulative discounted cost of the form

(20) $\qquad C^n = E \int_0^\infty e^{-\gamma t} \tilde{L}(\hat{\Phi}^n(t), \hat{\Psi}^n(t))\, dt,$

over all work-conserving admissible SCPs. Under the assumption that SCPs are work conserving, it is more convenient to work with the function $L : \mathbb{R}^k \times \mathbb{S}^k \to \mathbb{R}_+$ defined as

(21) $\qquad L(x, u) = \tilde{L}((\mathbb{1} \cdot x)^+ u,\ x - (\mathbb{1} \cdot x)^+ u).$

If work conservation holds, $(\mathbb{1} \cdot X^n - n)^+$ is equal to the number of customers waiting in all queues, namely $\mathbb{1} \cdot \Phi^n$. If $u^n \in \mathbb{S}^k$ denotes the proportion of customers of the different classes that are waiting in the queues, then

(22) $\qquad \Phi^n = (\mathbb{1} \cdot X^n - n)^+ u^n, \qquad \Psi^n = X^n - (\mathbb{1} \cdot X^n - n)^+ u^n.$

Hence (21) is merely a change of variables from $(\hat{\Phi}^n, \hat{\Psi}^n)$ to $(\hat{X}^n, u^n)$. The following will be assumed on $L$ and $\tilde{L}$.

ASSUMPTION 2. (i) $L(x, u) \geq 0$, $(x, u) \in \mathbb{R}^k \times \mathbb{S}^k$.

(ii) The mapping $(\phi, \psi) \mapsto \tilde{L}(\phi, \psi)$ is continuous. In particular, the mapping $(x, u) \mapsto L(x, u)$ is continuous.

(iii) There is $\varrho \in (0, 1)$ such that, for any compact $A \subset \mathbb{R}^k$,
$$|L(x, u) - L(y, u)| \leq c\|x - y\|^{\varrho}$$
holds for $u \in \mathbb{S}^k$ and $x, y \in A$, where $c$ depends only on $A$.

(iv) There are constants $c > 0$ and $m_L \geq 0$ such that $L(x, u) \leq c(1 + \|x\|^{m_L})$, $u \in \mathbb{S}^k$, $x \in \mathbb{R}^k$.

By applying an analogous change of variables to the state equations, both for the queueing system and for the diffusion, one can obtain these equations in a new form as follows. Equation (13) for $\hat{X}^n$ under work conservation takes the form

$$\hat{X}_t^n = \hat{X}^{0,n} + r_i \hat{W}_t^n + \int_0^t b^n(\hat{X}_s^n, u_s^n) \, ds, \tag{23}$$

where

$$b^n(\hat{X}, u) = \ell^n + (\mu^n - \theta^n)(\mathbb{1} \cdot \hat{X})^+ u - \mu^n \hat{X}, \tag{24}$$

and $r = \text{diag}(r_i; i \in K)$, $\ell^n = (\ell_1^n, \ldots, \ell_k^n)'$, $\mu^n = \text{diag}(\mu_i^n; i \in K)$, $\theta^n = \text{diag}(\theta_i^n; i \in K)$. Similarly, (16) for the diffusion model is now given as

$$X(t) = x + rW(t) + \int_0^t b(X(s), u(s)) \, ds, \tag{25}$$

where for $X \in \mathbb{R}^k$ and $u \in \mathbb{S}^k$,

$$b(X, u) = \ell + (\mu - \theta)(\mathbb{1} \cdot X)^+ u - \mu X, \tag{26}$$

and $\ell = (\ell_1, \ldots, \ell_k)'$, $\mu = \text{diag}(\mu_i; i \in K)$ and $\theta = \text{diag}(\theta_i; i \in K)$.

2.5. *Diffusion control problem.* Below we formulate a stochastic control problem for the minimization of

$$C = E \int_0^\infty e^{-\gamma t} L(X(t), u(t)) \, dt,$$

where $X$ is a controlled diffusion given by (25) over an appropriate class of control processes $u$, taking values in $\mathbb{S}^k$. We then state our first main result that there exists a measurable function $h: \mathbb{R}^k \to \mathbb{S}^k$ such that, upon setting $u_t = h(X_t)$, $t \geq 0$, the infimum in the problem is achieved.

DEFINITION 4. (i) We call
$$\pi = (\Omega, F, (F_t), P, u, W)$$
an *admissible system* if:

1. $(\Omega, F, (F_t), P)$ is a complete filtered probability space,
2. $u$ is a $\mathbb{S}^k$-valued, $F$-measurable, $(F_t)$-progressively measurable process, and $W$ is a standard $k$-dimensional $(F_t)$-Brownian motion.

The process $u$ is said to be a *control* associated with $\pi$.

(ii) We say that $X$ is a *controlled process* associated with initial data $x \in \mathbb{R}^k$ and an admissible system $\pi = (\Omega, F, (F_t), P, u, W)$, if:
1. $X$ is a continuous process on $(\Omega, F, P)$, $F$-measurable, $(F_t)$-adapted,
2. $\int_0^t |b(X(s), u(s))| \, ds < \infty$ for every $t \geq 0$, $P$-a.s. [recall that $b$ is defined in (26)],
3.
$$(27) \qquad X(t) = x + rW(t) + \int_0^t b(X(s), u(s)) \, ds, \qquad 0 \leq t < \infty,$$

holds $P$-a.s.

Proposition 2 shows that there is a unique controlled process $X$ associated with any $x$ and $\pi$. With an abuse of notation we sometimes denote the dependence on $x$ and $\pi$ by writing $P_x^\pi$ in place of $P$ and $E_x^\pi$ in place of $E$. Denote by $\Pi$ the class of all admissible systems.

PROPOSITION 2. *Let initial data $x \in \mathbb{R}^k$ and an admissible system $\pi \in \Pi$ be given. Then there exists a controlled process $X$ associated with $x$ and $\pi$. Moreover, if $X$ and $\bar{X}$ are controlled processes associated with $x$ and $\pi$, then $X(t) = \bar{X}(t)$, $t \geq 0$, $P$-a.s.*

For a proof see the Appendix.

For $x \in \mathbb{R}^k$ and $\pi \in \Pi$, let $X$ be the associated controlled process, and consider the cost function

$$C(x, \pi) = E_x^\pi \int_0^\infty e^{-\gamma t} L(X(t), u(t)) \, dt.$$

The value function $V$ for the control problem is defined as

$$V(x) = \inf_{\pi \in \Pi} C(x, \pi).$$

DEFINITION 5. Let $x \in \mathbb{R}^k$ be given. We say that a measurable function $h : \mathbb{R}^k \to \mathbb{S}^k$ is a *Markov control policy* if there is an admissible system $\pi$ and a controlled process $X$ corresponding to $x$ and $\pi$, such that $u_s = h(X_s)$, $s \geq 0$, $P$-a.s. We say that an admissible system $\pi$ is *optimal* for $x$, if $V(x) = C(x, \pi)$. We say that a Markov control policy is optimal for $x$ if the corresponding admissible system is.

The following constitutes a part of the first main result of this paper. Its full version that also characterizes the value function $V$ as the solution to an HJB equation, Theorem 3, is stated and proved in Section 3.

THEOREM 1. *Assume $L$ is continuous and satisfies Assumption 2(i), (iii) and (iv). Then there exists a Markov control policy, $h:\mathbb{R}^k \to \mathbb{S}^k$, which is optimal for all $x \in \mathbb{R}^k$.*

Throughout, $h$ denotes the function from Theorem 1.

2.6. *SCPs emerging from the diffusion control problem.* We formulate three SCPs that are based on the function $h$, and state our second main result, namely that these policies are, in an appropriate sense, asymptotically optimal.

*A P-SCP.* For each $t$, $\Psi^n(t)$ will be determined as a function of $X^n(t)$ only. Given $X^n(t)$, the diffusion control problem suggests setting

$$\Phi^n(t) = (\mathbb{1} \cdot X^n(t) - n)^+ h(\hat{X}^n(t)), \tag{28}$$

where as before

$$\hat{X}^n(t) = n^{1/2}\left(\frac{1}{n}X^n(t) - \rho\right).$$

There are two points, however, to which one must pay attention. First, the components of $\Phi^n$ must be integer-valued, in order to represent queue lengths; and second, the components of $\Psi^n = X^n - \Phi^n$ must be nonnegative, so that one serves only those customers present in the system.

For the first point, we need any measurable map $\Theta : \{y \in \mathbb{R}_+^k : \mathbb{1} \cdot y \in \mathbb{Z}\} \to \mathbb{Z}_+^k$ that preserves sums of components and introduces an error uniformly bounded by a constant, so that

$$\Phi^n(t) = \Theta[(\mathbb{1} \cdot X^n(t) - n)^+ h(\hat{X}^n(t))] \tag{29}$$

can be used in place of (28). For concreteness, take the following map. For $y \in \mathbb{R}_+^k$, write $y_i = \lfloor y_i \rfloor + \delta_i$, and set $z = \Theta(y)$ defined as $z_i = \lfloor y_i \rfloor$, $i = 1, \ldots, k-1$, and $z_k = y_k + \sum_{i=1}^{k-1} \delta_i$. Clearly, $\mathbb{1} \cdot z = \mathbb{1} \cdot y$, and $z \in \mathbb{Z}_+^k$ whenever $\mathbb{1} \cdot y \in \mathbb{Z}$. Moreover, $\|y - z\| \leq 2k$:

$$\|\Theta(y) - y\| \leq 2k, \qquad y \in \mathbb{R}_+^k. \tag{30}$$

For the second point, note that (29) might set $\Psi^n = X^n - \Phi^n$ in such a way that $\Psi^n$ is not in $\mathbb{R}_+^k$. For example, if $X^n = (n+1)e_1$ and $h(\hat{X}^n) = e_2$, then $\Phi_2^n = 1$, which means that $\Psi_2^n = -1$. Such a problem does not occur if

$$X_i^n(t) \geq (\mathbb{1} \cdot X^n(t) - n)^+ \qquad \forall i \in K. \tag{31}$$

When the problem does happen, the policy may be defined quite arbitrarily, subject only to being work conserving. For concreteness, when (31) is not met, we set $\Phi^n$ in accordance with a priority policy, where class $i$ receives priority $i$ (the higher $i$, the higher the priority). When (31) is met, we set $\Phi^n(t)$ as in (29). Finally, set $\Psi^n(t) = X^n(t) - \Phi^n(t)$, or equivalently, $\hat{\Psi}^n(t) = \hat{X}^n(t) - \hat{\Phi}^n(t)$. One verifies that the constraints (6) hold by construction.

We remark that the results of Section 4 will establish that (31) typically holds. This is basically due to the fact that the RHS, which represents the total number of customers waiting to be served, behaves at most as $O(n^{1/2})$, while the LHS, representing the number of customers at each class, is $O(n)$.

We next describe two alternative rules for determining sequences of N-SCPs.

*N-SCP* (i). To describe an N-SCP for each $n$, one needs to determine $\Psi^n$ so that the process $B^n(t)$ is nondecreasing. We describe a work-conserving SCP. A customer that arrives when there is a free server is instantaneously routed to a server. When a server becomes free, and there is at least one customer in the queue, we use the following scheme to determine which class to route to the server. This is in fact all that is to be determined. We look again at

$$M^n(t) := (\mathbb{1} \cdot X^n(t) - n)^+ h(\hat{X}^n(t)),$$

and consider the set $K^0$ of $i \in K$ for which $\Phi_i^n(t) \geq M_i^n(t) \vee 1$. Note that if there is at least one $i \in K$ with $\Phi_i^n \geq 1$, then $K^0$ is not empty. Indeed, suppose that $K^0$ is empty, and let $K' = \{i \in K : \Phi_i^n(t) \geq 1\}$. Then for $i \in K'$, $\Phi_i^n(t) < M_i^n(t)$. Hence by (18),

$$\mathbb{1} \cdot M^n = \mathbb{1} \cdot \Phi^n = \sum_{i \in K'} \Phi_i^n < \sum_{i \in K'} M_i^n \leq \mathbb{1} \cdot M^n,$$

a contradiction. We now choose the largest $i$ in $K^0$. Then a customer of class $i$ is routed to the free server. This procedure is performed instantaneously.

In heuristic terms, the scheme described above attempts to drive the system towards nearly achieving an equality of the form (28). This is done by sending to service customers of classes $i$ for which $\Phi_i^n \geq M_i^n$, thus obtaining approximate equality between the quantities $\Phi^n$ and $M^n$. A justification of this heuristic is a part of the proof of the result below.

*N-SCP* (ii). The N-SCP is defined precisely as the N-SCP (i), except that, for each $n$, the function $h$ is replaced by a function $h_n$, which may vary with $n$.

By defining the interarrival times $U_i^n(j)$ via $\breve{U}_i(j)$ [cf. (2)], we have assumed that they have finite variance. Here we strengthen this assumption.

ASSUMPTION 3. *Let $m_L$ be as in Assumption 2. Then there is a constant $m_U \geq 2$, $m_U > m_L$, such that $E(\check{U}_i(1))^{m_U} < \infty$.*

Our second main result is as follows.

THEOREM 2. *Let Assumptions 1–3 hold. Let $\hat{X}^{0,n} \in n^{-1/2}\mathbb{Z}^k$ be a sequence converging to $x \in \mathbb{R}^k$. Let a sequence of work-conserving admissible SCPs $\Psi^n$ be given, consider the corresponding processes $\Phi^n$, and let $\hat{\Phi}^n, \hat{\Psi}^n$ denote the corresponding rescaled processes.*

(i) *Let $\Psi^{n,*}, \Phi^{n,*}$ be a sequence as determined by the proposed P-SCP above, and let $\hat{\Phi}^{n,*}$ and $\hat{\Psi}^{n,*}$ be the corresponding rescaled processes. Then*

$$
\begin{aligned}
\lim_{n \to \infty} E \int_0^\infty e^{-\gamma t} \tilde{L}(\hat{\Phi}_t^{n,*}, \hat{\Psi}_t^{n,*}) \, dt \\
\leq \liminf_{n \to \infty} E \int_0^\infty e^{-\gamma t} \tilde{L}(\hat{\Phi}_t^n, \hat{\Psi}_t^n) \, dt.
\end{aligned}
\tag{32}
$$

*Moreover, the left-hand side is finite.*

(ii) *Assume that the restriction of $h$ to $\mathbb{X} := \{y \in \mathbb{R}^d : \mathbb{1} \cdot y > 0\}$ is locally Hölder continuous. Let $\Psi^{n,*}, \Phi^{n,*}$ be a sequence as determined by the proposed N-SCP (i) and $\hat{\Phi}^{n,*}, \hat{\Psi}^{n,*}$ be the corresponding rescaled processes. Then (32) holds.*

(iii) *Assume that the mapping $u \mapsto L(x, u)$ is convex on $\mathbb{S}^k$ for each $x \in \mathbb{R}^k$. Then there exists a sequence of functions $\{h_n\}$ with the following property. Let $\Psi^{n,*}, \Phi^{n,*}$ be a sequence as determined by the proposed N-SCP (ii), using the functions $\{h_n\}$, and $\hat{\Phi}^{n,*}, \hat{\Psi}^{n,*}$ be the corresponding rescaled processes. Then (32) holds.*

Item (i) of Theorem 2 establishes asymptotic optimality of the proposed sequence of preemptive SCPs, within all work-conserving SCPs. Item (ii) establishes asymptotic optimality of the proposed sequence of nonpreemptive SCPs, within all work-conserving SCPs, under the assumption that the function $h$ is locally Hölder continuous. In Proposition 3, we show that under some strict convexity assumptions on $L$, $h$ is locally Hölder continuous, and thus item (ii) applies. However, for linear costs, as $\tilde{L}(\hat{\Phi}, \hat{\Psi}) = c \cdot \hat{\Phi}$ ($c \in \mathbb{R}_+^k$ a constant), the resulting $h$ is discontinuous (see [21]), and this part of the theorem does not apply. Assuming only convexity of $L(x, \cdot)$, for each $x$ (which certainly holds for linear costs), item (iii) establishes asymptotic optimality of the proposed nonpreemptive SCPs, where $h$ is replaced by a sequence of functions $h_n$ that are locally Hölder continuous. Indeed, in Section 2.7 we discuss additional costs of interest, where $u \mapsto L(x, u)$ is convex for each $x$, implying that (i) and (iii) hold.

REMARK 2. The theorem is established by comparing both sides of (32) to the optimal cost in the corresponding diffusion control problem, denoted in Section 3 by $V(x)$. It is established below that, in fact, the left-hand side of (32) is equal to $V(x)$.

REMARK 3. As discussed in Section 5 (Corollary 1), for a sequence of N-SCPs that are not necessarily work conserving, Theorem 2 still holds given that work conservation is optimal among P-SCPs.

2.7. *Costs of interest.* The following result provides an example for a family of costs for which the assumptions on $h$ made in Theorem 2(ii) hold. It is proved in the Appendix.

PROPOSITION 3. *Let Assumption* 2 *hold, and assume that $\tilde{L}$ is of the form $\tilde{L}(\Phi, \Psi) = \sum_{i \in K} g_i(\Phi_i)$, where, for each $i \in K$, $g_i : [0, \infty) \to \mathbb{R}$ is in $C^2([0, \infty))$, and there is a constant $c_0 > 0$ such that $g_i'' \geq c_0$. Then the restriction of $h$ to $\mathbb{X}$ is locally Hölder continuous.*

Note that one can take in the above result $g_i(x) = c_i x^{p_i}$, $c_i > 0$, $p_i \geq 2$.

In the sequel we give examples of costs of interest, and specify the assumptions under which our main results apply. In all the cases below, $\tilde{L}$ and $L$ satisfy Assumption 2. Hence our results show asymptotic optimality of the proposed policies among *work-conserving* admissible policies.

*Queue lengths.* Let
$$\tilde{L}(\hat{\Phi}, \hat{\Psi}) = \ell(\hat{\Phi}),$$
where $\ell$ is nondecreasing as a function of $\hat{\Phi}_i$, for each $i$. It is assumed that $\ell \geq 0$ is locally Hölder continuous and satisfies a polynomial growth bound. Then
$$L(\hat{X}, u) = \ell((\mathbb{1} \cdot \hat{X})^+ u).$$

*Abandonment.* We need the following result, the proof of which is given in the Appendix.

LEMMA 1. *Under the assumptions of Theorem* 2, $ER_i^n(\mathring{T}_i^n(t)) = \theta_i E \mathring{T}_i^n(t)$.

The number of abandonments from queue $i$ up to time $t$, normalized by $\sqrt{n}$, is given by
$$\tilde{R}_i^n(t) := n^{-1/2} R_i^n(\mathring{T}_i^n(t)).$$

Consider the cost

$$C^n = \sum_i c_i E \int_0^\infty e^{-\gamma t} d\tilde{R}_i^n(t)$$

(the dependence of $C^n$ on the SCP is not indicated in this notation). Integrating by parts, using $e^{-\gamma t} E \tilde{R}_i^n(t) \to 0$ as $t \to 0$ and Lemma 1,

$$C^n = \sum_i \gamma c_i E \int_0^\infty e^{-\gamma t} \tilde{R}_i^n(t) \, dt$$

$$= \sum_i \gamma c_i \theta_i E \int_0^\infty e^{-\gamma t} \int_0^t \hat{\Phi}_i^n(s) \, ds \, dt$$

$$= E \int_0^\infty e^{-\gamma t} \left[ \sum_i c_i \theta_i \hat{\Phi}_i^n(t) \right] dt.$$

This is a special case of the queue-length cost considered in the previous paragraph.

*Delay.* For each of the customers $l$ ever present in the system, let $\mathrm{cl}(l)$ denote the class to which $l$ belongs, and let $\nu(l)$ denote the set of times at which customer $l$ is in the queue. We are interested in the cost

$$C^n = n^{-1/2} E \sum_l c_{\mathrm{cl}(l)} \int_{\nu(l)} e^{-\gamma t} \, dt,$$

where $c_i > 0$, $i \in K$, are constants. Since clearly, $\hat{\Phi}_i^n(t) = n^{-1/2} \sum \mathbf{1}_{t \in \nu(l)}$, where the sum extends over all class-$i$ customers $l$,

$$C^n = E \int_0^\infty e^{-\gamma t} \left[ \sum_{i \in K} c_i \hat{\Phi}_i^n(t) \right] dt.$$

This again can be treated within the framework of queue-length costs.

*Idling servers.* The number of servers that idle at time $t$ is given by $n - \mathbb{1} \cdot \Psi^n(t)$. With an appropriate normalization and discounting, this becomes

$$C^n = -E \int_0^\infty e^{-\gamma t} \mathbb{1} \cdot \hat{\Psi}^n(t) \, dt.$$

The corresponding costs are $\tilde{L}(\hat{\Phi}, \hat{\Psi}) = -\mathbb{1} \cdot \hat{\Psi}$ and $L(\hat{X}, u) = (\mathbb{1} \cdot \hat{X})^-$.

*Number of customers in the system.* The cost associated with the weighted normalized number of customers in the system is

$$C^n = E \int_0^\infty e^{-\gamma t} \sum_i c_i X_i(t) \, dt.$$

### 3. Stochastic control and the HJB equation.

3.1. *Moment estimates.* We begin with a key estimate for the results of this section.

PROPOSITION 4. *For any admissible system $\pi$, any $x, \bar{x} \in \mathbb{R}^k$, and corresponding controlled processes $X$ (associated with $x$ and $\pi$) and $\bar{X}$ (associated with $\bar{x}$ and $\pi$), the following hold:*

(i)
$$|X_t - \bar{X}_t| \leq |x - \bar{x}|(1 + e^{ct}), \qquad t \geq 0,$$

*$P$-a.s., where the constant $c$ does not depend on $\pi, x, \bar{x}$ and $t$.*

(ii) *For $m \in \mathbb{N}$,*
$$E_x^\pi |X(t)|^m \leq c_m (1 + \|x\|^m)(1 + t^m), \qquad t \geq 0,$$

*where the constants $c_1, c_2, \ldots$ do not depend on $\pi$, $x$ and $t$.*

PROOF. (i) Note that $|X(t) - \bar{X}(t)| \leq |x - \bar{x}| + c \int_0^t |X(s) - \bar{X}(s)| \, ds$, where $c$ is the Lipschitz constant for $x \mapsto b(x, u)$. The result follows from Gronwall's lemma.

(ii) Write $\Psi(t) = X - (\mathbb{1} \cdot X)^+ u$ and $\Psi_i(t) = \Psi(t) \cdot e_i$. Note that $\Psi_i(t) \leq X_i(t)$, and

(33) $$\sum_i \Psi_i(t) = 0 \wedge \sum_i X_i(t).$$

Then
$$X_i(t) = x_i + r_i W_i(t) + \int_0^t [-\theta_i X_i(s) - (\mu_i - \theta_i)\Psi_i(s) + \ell_i] \, ds, \qquad i \in K, \ t \geq 0.$$

Let $K_1$ be the set of $i \in K$, where $\mu_i \geq \theta_i$, and $K_2 = K \setminus K_1$. Define, for each $i$, $\Upsilon_i$ as the unique solution (cf. Theorems 5.2.5 and 5.2.9 of [25]) to the equation
$$\Upsilon_i(t) = x_i + r_i W_i(t) + \int_0^t [-\mu_i \Upsilon_i(s) + \ell_i] \, ds.$$

Then $X_i - \Upsilon_i$ is differentiable, $X_i(0) - \Upsilon_i(0) = 0$, and for $i \in K_1$,
$$\frac{d}{dt}(X_i(t) - \Upsilon_i(t)) = -\theta_i X_i - (\mu_i - \theta_i)\Psi_i + \mu_i \Upsilon_i$$
$$\geq -\mu_i(X_i(t) - \Upsilon_i(t)).$$

Similarly, the reverse inequality holds when $i \in K_2$. By comparison of ODEs (Theorem I.7 in [5]),

(34) $X_i(t) \geq \Upsilon_i(t), \qquad i \in K_1; \qquad X_i(t) \leq \Upsilon_i(t), \qquad i \in K_2; \ t \geq 0 \quad$ a.s.

If $z$ is a vector satisfying the bounds $z_i \geq a_i$ for all $i$, and $\sum_i z_i \leq A$, then its norm can be bounded as follows:

$$\|z\| \leq \sum_i (z_i - a_i) + \|a\| \leq A + 2\|a\|. \tag{35}$$

We have in (34) inequalities analogous to $z_i \geq a_i$, when we consider $z_i = c_i X_i$, where $c_i > 0$, $i \in K_1$, and $c_i < 0$, $i \in K_2$. Below, we obtain an inequality analogous to $\sum_i z_i \leq A$, by finding an upper bound on the quantity $\sum_K c_i X_i$. To this end, note first that by (33) and (34),

$$\sum_{K_1} -\Psi_i + \sum_{K_2} (X_i - \Psi_i) = -\left(0 \wedge \sum_K X_i\right) + \sum_{K_2} X_i$$

$$= \mathbf{1}_{\{\sum_K X_i \geq 0\}} \sum_{K_2} X_i - \mathbf{1}_{\{\sum_K X_i < 0\}} \sum_{K_1} X_i \tag{36}$$

$$\leq \sum_K |\Upsilon_i|.$$

Next, let $c > 0$ be so small that $1 + c(1 - \theta_i/\mu_i) \geq 0$ for all $i \in K_2$. Then also $[1 + c(1 - \theta_i/\mu_i)](X_i - \Psi_i) \geq 0$, and as a result,

$$c\left[-\frac{\theta_i}{\mu_i}(\Psi_i - X_i) + \Psi_i\right] \leq (X_i - \Psi_i) + cX_i, \qquad i \in K_2. \tag{37}$$

Hence, denoting $\tilde{x} = \sum_{K_1} \mu_i^{-1} x_i - \sum_{K_2} c\mu_i^{-1} x_i$, $\tilde{W}(t) = \sum_{K_1} r_i W_i(t) - \sum_{K_2} cr_i W_i(t)$ and $\tilde{\ell} = \sum_{K_1} \ell_i - \sum_{K_2} c\ell_i$, we have by (34), (36) and (37)

$$\sum_{K_1} \mu_i^{-1} X_i(t) - \sum_{K_2} c\mu_i^{-1} X_i(t)$$

$$= \tilde{x} + \tilde{W}(t) + \tilde{\ell}t + \int_0^t \sum_{K_1} [\mu_i^{-1}\theta_i(\Psi_i(s) - X_i(s)) - \Psi_i(s)]\, ds$$

$$+ \int_0^t \sum_{K_2} c[-\mu_i^{-1}\theta_i(\Psi_i(s) - X_i(s)) + \Psi_i(s)]\, ds$$

$$\leq \tilde{x} + \tilde{W}(t) + \tilde{\ell}t + \int_0^t \left[\sum_{K_1} -\Psi_i(s) + \sum_{K_2} [(X_i(s) - \Psi_i(s)) + cX_i(s)]\right] ds$$

$$\leq \tilde{x} + \tilde{W}(t) + \tilde{\ell}t + \int_0^t \left[\sum_K |\Upsilon_i(s)| + \sum_{K_2} c|\Upsilon_i(s)|\right] ds.$$

Denoting $Z_i = c_i X_i$, where $c_i = \mu_i^{-1}$, $i \in K_1$, and $c_i = -c\mu_i^{-1}$, $i \in K_2$, we have from (34), (35) and the above, that, for some positive constants $C_1, C_2$,

$$C_1 \|X(t)\| \leq \|Z(t)\| \leq \tilde{x} + \tilde{W}(t) + \tilde{\ell}t + C_2 \int_0^t \sum_K |\Upsilon_i(s)|\, ds + C_2 \sum_K |\Upsilon_i(t)|.$$

It is easy to show that $E|\Upsilon_i(t)|^2 \leq C(1+|x_i|^2)$, for some constant $C$, and since $\Upsilon_i$ are Gaussian, $E|\Upsilon_i(t)|^m \leq \tilde{c}_m(1+|x_i|^m)$ for $m=1,2,\ldots$. It easily follows that

$$(38) \quad E\|X(t)\|^m \leq c_m(1+\|x\|^m)(1+t^m). \qquad \square$$

REMARK 4. We record a consequence of the proof to be used in Section 4. Recall (23) which holds under the work-conservation condition (19). Arguing analogously to the proof of Proposition 4, under (19) one obtains

$$(38) \quad \|\hat{X}^n(t)\| \leq c\bigg[\|\hat{X}^{0,n}\| + \|\hat{W}_t^n\| + t + \int_0^t \|\Upsilon^n(s)\|\,ds + \|\Upsilon^n(t)\|\bigg],$$

where $c$ does not depend on $n$ or $t$, and where $\Upsilon^n$ is the unique solution to

$$(39) \quad \Upsilon_i^n(t) = \hat{X}_i^{0,n} + r_i^n \hat{W}_i^n(t) + \int_0^t (-\mu_i^n \Upsilon_i^n(s) + \ell_i^n)\,ds.$$

3.2. *Cost and value.* Recall that for $x \in \mathbb{R}^k$ and $\pi \in \Pi$, the cost and value are defined as

$$C(x,\pi) = E_x^\pi \int_0^\infty e^{-\gamma t} L(X(t), u(t))\,dt,$$

$$V(x) = \inf_{\pi \in \Pi} C(x,\pi).$$

We assume in this section that $L(x,u)$ satisfies Assumption 2, except that part (ii) should be understood as the assumption that $L$ is continuous (the notation $\tilde{L}$ is not needed in this section).

To state the next result, we need to formulate a control problem on a bounded domain. In the sequel, $\Gamma$ will denote a bounded open connected subset of $\mathbb{R}^k$ with smooth (say, $C^\infty$) boundary. Let $g: \mathbb{R}_+ \times \partial\Gamma \to \mathbb{R}_+$ be a continuous function. For $x \in \Gamma$ and $\pi \in \Pi$, we define

$$C_{\Gamma,g}(x,\pi) = E_x^\pi\bigg[\int_0^\tau e^{-\gamma t} L(X_t, u_t)\,dt + g(\tau, X_\tau)\bigg],$$

where $X$ is the corresponding controlled process, and

$$\tau = \inf\{t : X_t \notin \Gamma\}.$$

We also let

$$V_{\Gamma,g}(x) = \inf_{\pi \in \Pi} C_{\Gamma,g}(x,\pi).$$

PROPOSITION 5. *Assume $L$ is continuous and satisfies Assumption* 2(i), (iii) *and* (iv). *Then:*

(i) *There is a constant $c$ such that $V(x) \leq c(1+\|x\|^{m_L})$, $x \in \mathbb{R}^k$.*
(ii) *$V$ is continuous on $\mathbb{R}^k$.*

(iii) *Let $\Gamma \subset \mathbb{R}^k$ be a smooth domain. Let $g(t,x) = e^{-\gamma t}V(x)$ for $t \geq 0$ and $x \in \partial\Gamma$. Then $V = V_{\Gamma,g}$ in $\Gamma$.*

PROOF. (i) This is immediate from the polynomial growth condition on $L$ and Proposition 4(ii).

(ii) Fix an arbitrary open ball of radius 1, $\nu = B(y,1)$. Let $x \in \nu$ be given, and for $\varepsilon > 0$, let $\pi = (\Omega, F, (F_t), P, u, W)$ be such that

$$C(x,\pi) \leq V(x) + \varepsilon.$$

Let $X$ be the controlled process associated with $x$ and $\pi$. Let $\bar{X}$ be the controlled process on the same probability space, associated with $\pi$ and some $\bar{x} \in \nu$. Denote $m = m_L$ (as in Assumption 2). Let $A(T) = B(y, T^{2m+3})$. Let $c_1(T)$ be the Hölder constant for $L$ on $A(T)$. By Proposition 4(ii),

$$E_z^\pi \|X(t)\|^m \leq \hat{c}(1+t^m), \qquad t \geq 0, \ z \in \nu, \tag{40}$$

where $\hat{c} = \hat{c}(\nu)$. Then for any $T \in [1,\infty)$ and $t \in [0,T]$, Proposition 4(i), (40) and the Cauchy–Schwarz inequality imply

$$\begin{aligned}
E|L(X_t, u_t) &- L(\bar{X}_t, u_t)| \\
&\leq c_1(T) E[\mathbf{1}_{\{X_t, \bar{X}_t \in A(T)\}} \|X_t - \bar{X}_t\|^\varrho] \\
&\quad + cE[\mathbf{1}_{\{\text{either } X_t \text{ or } \bar{X}_t \notin A(T)\}}(1 + \|X_t\|^m + \|\bar{X}_t\|^m)] \\
&\leq c_1(T)(1 + e^{cT})\|x - \bar{x}\|^\varrho + c[p(T) + \bar{p}(T)]^{1/2}\hat{c}(1 + T^m),
\end{aligned}$$

where

$$p(T) = \sup_{s \leq T} P(X_s \notin A(T)), \qquad \bar{p}(T) = \sup_{s \leq T} P(\bar{X}_s \notin A(T)).$$

The moment bounds on $\|X_t\|$ imply that

$$p(T) + \bar{p}(T) \leq c_2(\nu) T^{-2m-2},$$

where $c_2(\nu)$ depends on $\nu$, but not on $x, \bar{x} \in \nu$. Hence, writing $c_3(\nu) = \hat{c}c_2(\nu)^{1/2}$,

$$\begin{aligned}
C(x,\pi) &- C(\bar{x},\pi) \\
&= E\int_0^\infty e^{-\gamma t}(L(X_t, u_t) - L(\bar{X}_t, u_t))\,dt \\
&\leq \{c_1(T)(1 + e^{cT})\|x - \bar{x}\|^\varrho + cc_3(\nu)T^{-m-1}(1 + T^m)\} E\int_0^T e^{-\gamma t}\,dt \\
&\quad + c\int_T^\infty e^{-\gamma t}(1 + E\|X_t\|^m + E\|\bar{X}_t\|^m)\,dt \\
&\leq c_4(T)\|x - \bar{x}\|^\varrho + c_5(\nu)\alpha(T),
\end{aligned}$$

where $c_4(T)$ depends only on $T$, $c_5(\nu)$ depends only on $\nu$, and $\alpha(T) \to 0$ as $T \to \infty$. Let $T$ be so large that $c_5(\nu)\alpha(T) \leq \varepsilon$. Next choose $\delta > 0$ so small that $\{0 < \|x - \bar{x}\| < \delta$ and $x, \bar{x} \in \nu\}$ implies $c_4(T)\|x - \bar{x}\|^\varrho \leq \varepsilon$. Then for such $x, \bar{x}$ one has $V(\bar{x}) \leq C(\bar{x}, \pi) \leq C(x, \pi) + 2\varepsilon \leq V(x) + 3\varepsilon$. Note that the choice of $\delta$ does not depend on $x, \bar{x}$ (in particular, it does not depend on $\pi$!). Therefore, the inequality $V(\bar{x}) \leq V(x) + 3\varepsilon$ holds for all $x, \bar{x} \in \nu$ for which $\|x - \bar{x}\| < \delta$. This shows that $V$ is continuous.

(iii) This is a standard result (the principle of optimality), which, in the current context, can be proved similarly to the results of [6], Section III.1. □

3.3. *The HJB equation and optimality.* The HJB equation associated with the stochastic control problem is (cf. [11])

$$\mathcal{L}f + H(x, Df) - \gamma f = 0, \tag{41}$$

where $\mathcal{L} = (1/2)\sum_i r_i^2 \partial^2/\partial x_i^2$, and

$$H(x, p) = \inf_{u \in \mathbb{S}^k}[b(x, u) \cdot p + L(x, u)].$$

The equation is considered on $\mathbb{R}^k$ with the growth condition

$$\exists C, m, \quad |f(x)| \leq C(1 + \|x\|^m), \quad x \in \mathbb{R}^k. \tag{42}$$

We say that $f$ is a solution to (41) if it is of class $C^2$, and the equation is satisfied everywhere in $\mathbb{R}^k$.

THEOREM 3. *Assume $L$ is continuous and satisfies Assumption* 2(i), (iii) *and* (iv). *Then there exists a classical solution $f \in C^{2,\varrho}_{\mathrm{pol}}(\mathbb{R}^k)$ to* (41), (42), *and this solution is unique in $C^2_{\mathrm{pol}}(\mathbb{R}^k)$. Moreover, the value $V$ is equal to $f$. Furthermore, there exists a Markov control policy which is optimal for all $x \in \mathbb{R}^k$.*

PROOF. We first consider equation (41) on a smooth open bounded connected domain $\Gamma$, satisfying an exterior sphere condition, with boundary conditions

$$f(x) = V(x), \quad x \in \partial\Gamma. \tag{43}$$

The key is a result from [16] regarding existence of classical solutions in bounded domains, with merely continuous boundary conditions. To use this result, we verify the following two conditions:

(i) $|H(x, p)| \leq c(1 + \|p\|)$ for $x \in \Gamma$, where $c$ does not depend on $x$ or $p$.
(ii) $H(x, p) \in C^\varepsilon(\overline{\Gamma} \times \mathbb{R}^k)$, some $\varepsilon \in (0, 1)$.

Item (i) is immediate from the local boundedness of $b(x,u)$ and $L(x,u)$. Next we show that item (ii) holds. For $\delta > 0$, let $v$ be such that $H(y,q) \geq b(y,v) \cdot q + L(y,v) - \delta$. Write

$$H(x,p) - H(y,q) \leq b(x,v) \cdot p + L(x,v) - b(y,v) \cdot q - L(y,v) + \delta.$$

Using the Hölder property of $L$ in $x$ uniformly for $(x,v) \in \overline{\Gamma} \times \mathbb{S}^k$, and the Lipschitz property of $b$ in $x$, uniformly in $(x,v)$,

$$H(x,p) - H(y,q) \leq c\|p - q\| + c\|p\|\|x - y\| + c\|x - y\|^\varrho + \delta.$$

Since $\delta > 0$ is arbitrary, it can be dropped. This shows that $H$ is Hölder continuous with exponent $\varrho$, uniformly over compact subsets of $\overline{\Gamma} \times \mathbb{R}^k$. Hence (ii) holds.

Defining for $(x,z,p) \in \Gamma \times \mathbb{R} \times \mathbb{R}^k$, $A(x,z,p) = (1/2)r^2 p$, $B(x,z,p) = H(x,p) - \gamma z$, one can write (41) in divergence form as

$$\text{div}\, A(x,f,Df) + B(x,f,Df) = 0.$$

The hypotheses of Theorem 15.19 of [16] regarding the coefficients $A$ and $B$ hold in view of (i) and (ii). Indeed, $B$ is Hölder continuous of exponent $\varrho$, uniformly on compact subsets of $\Gamma \times \mathbb{R} \times \mathbb{R}^k$. Moreover, with $\tau = 0$, $\nu(z) = (1/2)\min_i r_i^2$, $\mu(z) = c(1 + \|z\|)$, $\alpha = 2$, $b_1 = 0$ and $a_1 = 0$, one checks that the conditions (15.59), (15.64), (15.66) and (10.23) of [16] are satisfied. Theorem 15.19 of [16] therefore applies. [We comment that there is a typo in the statement of the conditions of the theorem in [16]: the reference should be to condition (15.59) instead of (15.60).] It states that there exists a solution to (41) in $C^{2,\varrho}(\Gamma) \cap C(\overline{\Gamma})$, satisfying the continuous boundary condition (43). We denote this solution by $f$.

Let $x \in \Gamma$. Let $\pi$ be any admissible system and let $X$ be the controlled process associated with $x$ and $\pi$. Let $\tau$ denote the first time $X$ hits $\partial \Gamma$. Using Itô's formula for the $C^{1,2}(\mathbb{R}_+ \times \Gamma)$ function $e^{-\gamma t} f(x)$, in conjunction with the inequality

$$\mathcal{L}f(y) + b(y,u) \cdot Df(y) + L(y,u) - \gamma f(y) \geq 0, \qquad y \in \Gamma, \ u \in \mathbb{S}^k,$$

satisfied by $f$, one obtains

$$(44) \quad \begin{aligned} f(x) &\leq \int_0^{t\wedge\tau} e^{-\gamma s} L(X_s, u_s)\, ds \\ &\quad + e^{-\gamma(t\wedge\tau)} f(X_{t\wedge\tau}) - \int_0^{t\wedge\tau} e^{-\gamma s} Df(X_s) \cdot r\, dW_s. \end{aligned}$$

Taking expectation and then sending $t \to \infty$, using the monotone convergence theorem as well as the bounded convergence theorem, we have with $g(t,x) = e^{-\gamma t} V(x)$,

$$f(x) \leq E_x^\pi \left[ \int_0^\tau e^{-\gamma s} L(X_s, u_s)\, ds + e^{-\gamma \tau} V(X_\tau) \right] = C_{\Gamma,g}(x,\pi).$$

Taking the infimum over $\pi \in \Pi$, we have

$$f(x) \leq V_{\Gamma,g}(x) = V(x), \qquad x \in \Gamma,$$

where the last equality follows from Proposition 5(iii).

In order to obtain the equality $f = V$ on $\Gamma$, we next show there exist optimal Markov control policies for the control problem on $\Gamma$. Let

(45) $$\varphi(x,u) = b(x,u) \cdot Df(x) + L(x,u), \qquad x \in \Gamma, u \in \mathbb{S}^k.$$

Note that $\varphi$ is continuous on $\Gamma \times \mathbb{S}^k$. For each $x$, consider the set $U_x \neq \varnothing$ of $u \in \mathbb{S}^k$ for which

$$\varphi(x,u) = \inf_{v \in \mathbb{S}^k} \varphi(x,v).$$

We show that there exists a measurable selection of $U_x$, namely there is a measurable function $h$ from $(\Gamma, \mathcal{B}(\Gamma))$ to $(\mathbb{S}^k, \mathcal{B}(\mathbb{S}^k))$ with $h(x) \in U_x$, $x \in \Gamma$.

Let $x_n \in \Gamma$ and assume $\lim_n x_n = x \in \Gamma$. Let $u_n$ be any sequence such that $u_n \in U_{x_n}$. We claim that any accumulation point of $u_n$ is in $U_x$, for if this is not true, then by continuity of $\varphi$, there is a converging subsequence $u_m$, converging to $\bar{u}$, and there is a $\hat{u}$ such that $\delta := \varphi(x,\bar{u}) - \varphi(x,\hat{u}) > 0$. Hence for all $m$ large, $\varphi(x_m, u_m) \geq \varphi(x, \hat{u}) + \delta/2 \geq \varphi(x_m, \hat{u}) + \delta/4$, contradicting $u_m \in U_{x_m}$.

As a consequence, the assumptions of Corollary 10.3 in the Appendix of [10] are satisfied, and it follows that there exists a measurable selection $h: \Gamma \to \mathbb{S}^k$ of $(U_x, x \in \Gamma)$.

We extend $h$ to $\mathbb{R}^k$ in a measurable way so that it takes values in $\mathbb{S}^k$ (but is otherwise arbitrary). Clearly, $x \mapsto b(x, h(x))$ is measurable. Consider the autonomous SDE

(46) $$X(t) = x + rW(t) + \int_0^t \hat{b}(X_s)\,ds,$$

where $\hat{b}(y)$ agrees with $b(y, h(y))$ on $\Gamma$, and is set to zero off $\Gamma$. Then $\hat{b}$ is measurable and bounded on $\mathbb{R}^k$. By Proposition 5.3.6 of [25], there exists a weak solution to this equation. That is, there exists a complete filtered probability space on which $X$ is adapted and $W$ is a $k$-dimensional Brownian motion, such that (46) holds for $t \geq 0$, a.s. On this probability space, consider the process $u_s = h(X_s)$. Since $X$ has continuous paths and is adapted, it is progressively measurable (see Proposition 1.13 of [25]) and by measurability of $h$, so is $u$. Denote by $\pi$ the admissible system thus constructed. Then for $s < \tau$, $u_s \in U_{X_s}$ and

$$b(X_s, u_s) \cdot Df(X_s) + L(X_s, u_s) = H(X_s, Df(X_s)).$$

Hence

$$\mathcal{L}f(X) + b(X_s, u_s) \cdot Df(X_s) + L(X_s, u_s) - \gamma f(X_s) = 0, \qquad s < \tau.$$

A use of Itô's formula and the convergence theorems just as before now shows that

$$f(x) = E_x^\pi \left[ \int_0^\tau e^{-\gamma s} L(X_s, u_s)\, ds + e^{-\gamma \tau} V(X_\tau) \right] = C_{\Gamma,g}(x, \pi), \qquad x \in \Gamma,$$

with $g$ as above. This, together with the principle of optimality, shows that $f \geq V_{\Gamma,g} = V$ on $\Gamma$. Summarizing, $f = V$ on $\Gamma$.

In particular, $V \in C^{2,\varrho}(\Gamma)$ and is a classical solution to the HJB equation. $\Gamma$ can now be taken arbitrarily large, and this shows that $V \in C^{2,\varrho}(\mathbb{R}^k)$, and that it satisfies the HJB equation on $\mathbb{R}^k$. In view of Proposition 5(i), it also satisfies the polynomial growth condition. As a result, there exists a classical solution to (41) in $C^{2,\varrho}(\mathbb{R}^k)$, again denoted by $f$, satisfying (42), and moreover, $V = f$.

It remains to show uniqueness within $C^2_{\text{pol}}(\mathbb{R}^k)$ and existence of optimal Markov control policies for the problem on $\mathbb{R}^k$. Let $\bar{f} \in C^2_{\text{pol}}(\mathbb{R}^k)$ be any solution to (41). Then analogously to (44), we obtain

$$\bar{f}(x) \leq \int_0^t e^{-\gamma s} L(X_s, u_s)\, ds + e^{-\gamma t} \bar{f}(X_t) - \int_0^t e^{-\gamma s} D\bar{f}(X_s) \cdot r\, dW_s.$$

Taking expectation, sending $t \to \infty$, using the polynomial growth of $\bar{f}$ and the moment bounds on $\|X_t\|$, one has that $\bar{f}(x) \leq C(x, \pi)$, where $\pi \in \Pi$ is arbitrary. Consequently, $\bar{f} \leq V$ on $\mathbb{R}^d$.

The proof of existence of optimal Markov policies as well as the inequality $V \leq \bar{f}$ on $\mathbb{R}^k$ is completely analogous to that on $\Gamma$, where one replaces $\Gamma$ by $\mathbb{R}^k$ and uses again the polynomial growth condition of $\bar{f}$. The weak existence of solutions to (46) follows on noting that $\hat{b}$ satisfies a linear growth condition of the form $\|\hat{b}(y)\| \leq x(1 + \|y\|)$, $y \in \mathbb{R}^k$, and using again Proposition 5.3.6 of [25]. Hence $V = \bar{f}$ on $\mathbb{R}^k$. We conclude that $f$ is the unique solution in $C^2_{\text{pol}}(\mathbb{R}^k)$, that $V = f$, and that there exists a Markov control policy, optimal for all $x \in \mathbb{R}^k$. □

**4. Asymptotic optimality.** In this section we prove asymptotic optimality of the proposed SCPs. As in the statement of Theorem 2, all SCPs are assumed to be work conserving in this section. Recall from Section 2 that the processes $\Phi^n$ and $\Psi^n$ represent the number of customers waiting in each queue, and, respectively, the number of servers working on jobs of each class. Let $u^n$ be an $\mathbb{S}^k$-valued process, determined as

$$(47) \qquad u^n = \begin{cases} \Phi^n/(\mathbb{1} \cdot X^n - n)^+, & \mathbb{1} \cdot X^n - n > 0, \\ u_0, & \mathbb{1} \cdot X^n - n \leq 0, \end{cases}$$

where $u_0$ is some fixed, arbitrary element of $\mathbb{S}^k$. As in the paragraph preceding Assumption 2, $u^n$ represents the fraction of customers of each class that are waiting in the queues (whenever there are such customers). As a result one can determine $\hat{\Phi}^n$ and $\hat{\Psi}^n$ from $u^n$ and $\hat{X}^n$ as $\hat{\Phi}^n = (\mathbb{1} \cdot \hat{X}^n)^+ u^n$ and $\hat{\Psi}^n = \hat{X}^n - \hat{\Phi}^n$.

Throughout this section let $f$ denote the unique $C^2_{\text{pol}}$ solution to (41) (cf. Theorem 3). Let

$$(48) \quad K^n_t = b(\hat{X}^n_t, u^n_t) \cdot Df(\hat{X}^n_t) + L(\hat{X}^n_t, u^n_t) - H(\hat{X}^n_t, Df(\hat{X}^n_t)) \geq 0.$$

A condition that plays a central role in the convergence proof is

$$(49) \qquad \int_0^\cdot e^{-\gamma s} K^n_s \, ds \Rightarrow 0.$$

THEOREM 4. (i) *Let Assumptions* 1–3 *hold. Let* $\hat{X}^{0,n} \in n^{-1/2}\mathbb{Z}^k$ *be a sequence converging to* $x \in \mathbb{R}^k$. *Let a sequence of work-conserving admissible SCPs be given [namely,* (19) *holds], let* $\hat{X}^n$ *be the corresponding normalized controlled processes starting from* $\hat{X}^{0,n}$ *and let* $u^n$ *be given by* (47). *Then*

$$\liminf_{n \to \infty} E \int_0^\infty e^{-\gamma t} L(\hat{X}^n_t, u^n_t) \, dt \geq V(x).$$

(ii) *Assume, in addition, that* (49) *is satisfied. Then*

$$\limsup_{n \to \infty} E \int_0^\infty e^{-\gamma t} L(\hat{X}^n_t, u^n_t) \, dt \leq V(x).$$

In what follows we prove Theorem 4. We treat both parts (i) and (ii) simultaneously. Whenever there is a reference to part (ii), we indicate explicitly that (49) holds. It will be convenient to work with both representations (13) and (23) for $\hat{X}^n$ in this section. Denote

$$(50) \qquad Y^n_t = \int_0^t b^n(\hat{X}^n_s, u^n_s) \, ds, \qquad Z^n_t = \int_0^t e^{-\gamma s} L(\hat{X}^n_s, u^n_s) \, ds.$$

Let $(\mathcal{F}^n_t)$ be the filtration (9). Note that, by definition, the processes $\hat{X}^n, \hat{\Phi}^n, \hat{\Psi}^n$ are adapted to $\mathcal{F}^n$. Hence by (47) and (50), so are the processes $u^n, Y^n$ and $Z^n$.

LEMMA 2. *Under Assumption* 3,

$$E(\|\hat{A}^n\|^*_t)^{m_U} \leq c(1 + t^{m_U/2}), \qquad n \in \mathbb{N}, \ t \in \mathbb{R}_+,$$

*where $c$ does not depend on $n$ or $t$.*

PROOF. This is a consequence of Theorem 4 of [26], which, under the assumption $E(\check{U}_i(1))^{m_U} < \infty$, $m_U \geq 2$, states that

$$(51) \qquad E \sup_{s \leq t} |n^{-1/2}(A^n_i(ns) - n\lambda_i s)|^{m_U} \leq c(1 + t^{m_U/2}),$$

where

$$A_i(t) = \sup\left\{m \geq 0 : \sum_{j=1}^m \check{U}_i(j) \leq t\right\}, \qquad t \geq 0,$$

and $c$ does not depend on $n$ or $t$. Indeed, by (2) and (3), $A_i^n(t) = A_i(\lambda_i^n t)$. Let $C = \sup_n[\lambda_i^n/(n\lambda_i)]$ and note that $C < \infty$ by Assumption 1. Then

$$|\hat{A}_i^n|_t^* = \sup_{s \leq t} n^{-1/2} |A_i(\lambda_i^n s) - \lambda_i^n s|$$

$$\leq \sup_{s \leq Ct} n^{-1/2} |A_i(n\lambda_i s) - n\lambda_i s|.$$

The lemma follows from (51). □

LEMMA 3. *Under the assumptions of Theorem* 4(i), *the processes* $\hat{X}^n$ *satisfy* $E\|\hat{X}^n(t)\|^{m_U} \leq c(1 + \|x\|^{\bar{m}})(1 + t^{\bar{m}})$, *where* $\bar{m}$ *and* $c$ *do not depend on* $n$, $x$ *or* $t$.

PROOF. Since we are assuming work conservation, (38) applies. Solving for $\Upsilon^n$ of (39), we obtain

$$\Upsilon_i^n(t) = \hat{X}_i^{0,n} e^{-\gamma t} + \tilde{W}_i^n(t) - \mu_i^n \int_0^t \tilde{W}_i^n(s) e^{-\mu_i^n(t-s)} ds,$$

where

$$\tilde{W}_i^n(t) = r_i \hat{W}_i^n(t) + \ell_i^n t.$$

Hence

$$(52) \quad \|\hat{X}_t^n\| \leq c\left[1 + t^2 + \|\hat{X}^{0,n}\| + \|\hat{W}_t^n\| + \int_0^t \|\hat{W}_s^n\| ds + \int_0^t \int_0^s \|\hat{W}_\theta^n\| d\theta ds\right].$$

By (14), using $\bar{\Psi}_i^n(s) \leq 1$ and $\bar{\Phi}_i^n(s) \leq \xi^n(s) := \max_i[n^{-1}X_i^{0,n} + n^{-1}A_i^n(s)]$,

$$(53) \qquad \|\hat{W}^n(t)\| \leq \|\hat{A}^n(t)\| + \sup_{s \leq t} \|\hat{S}^n(s)\| + \sup_{s \leq \xi^n(t)} \|\hat{R}^n(s)\|.$$

Denote $p = m_U$. Apply Burkholder's inequality (cf. [33], page 175) to the (discontinuous) martingale $\hat{S}^n$, denoting by $[M]$ the quadratic variation processes associated with $M$, and recalling that if a process $M$ taking real values has sample paths of bounded variation, then $[M](t) = M_0^2 + \sum_{0 < s \leq t} (\Delta M_s)^2$. Denoting by $\chi_i^n(t)$ a Poisson random variable with parameter $n\mu_i^n t$ and using the convergence $\mu_i^n \to \mu_i$, we obtain

$$E\sup_{s \leq t} |\hat{S}_i^n(s)|^p \leq cE([\hat{S}_i^n](t))^{p/2}$$

$$= cn^{-p/2} E(\chi_i^n(t))^{p/2}$$

$$\leq c_p n^{-p/2}(n\mu_i^n t)^{p/2}$$
$$\leq ct^{p/2},$$

where $c$ does not depend on $n$ or $t$. Similarly, $E\sup_{s\leq t}|\tilde{R}_i^n(s)|^p \leq ct^{p/2}$. Therefore, by the independence of $A^n$ and $R^n$ and Assumption 3,

$$E\sup_{s\leq \xi^n(t)}|\hat{R}_i^n(s)|^p = E\left\{E\left[\sup_{s\leq \xi^n(t)}|\hat{R}_i^n(s)|^p\Big|\xi^n(t)\right]\right\}$$
$$\leq cE(\xi^n(t))^{p/2}$$
$$\leq c(1+t^q),$$

where $q$ does not depend on $n$ or $t$. Lemma 2 and an application of Minkowski's inequality to (53) show that there is $m$ not depending on $n$ or $t$ such that

(54) $$E\|\hat{W}^n(t)\|^p \leq c(1+t^m), \qquad t\geq 0.$$

The lemma now follows from (52). □

LEMMA 4. *Let the assumptions of Theorem 4(i) hold.*

(i) $(\hat{A}^n, \hat{S}^n, \hat{R}^n) \Rightarrow (A, S, R)$, *where $A, S$ and $R$ are independent Brownian motions with zero drift and variance matrices* $\mathrm{diag}(\lambda_i C_{U,i}^2)_{i\in K}$, $\mathrm{diag}(\mu_i)_{i\in K}$, *and, respectively,* $\mathrm{diag}(\theta_i)_{i\in K}$.

(ii) *One has*

(55) $$(\bar{\Psi}^n, \bar{\Phi}^n) \Rightarrow (\rho, 0) \qquad in \ (\mathbb{D}(\mathbb{R}^k))^2$$

[*the process that is constantly $(\rho, 0)$*].

(iii) *The sequence $(\hat{X}^n, Y^n, Z^n, \hat{W}^n)$ is tight* [*in $(\mathbb{D}(\mathbb{R}^k))^4$*].

PROOF. (i) By the assumption on the finite second moment and i.i.d. structure of the interarrival times, and by Assumption 1, the results of [22] imply (i).

(ii) Since work conservation (19) is assumed, we can use (38). Note that $\bar{X}^n - \rho = n^{-1/2}\hat{X}^n$. By part (i), $n^{-1/2}\hat{W}^n \Rightarrow 0$. Also, $n^{-1/2}\hat{X}^{0,n} \to 0$. Hence by Gronwall's lemma, $n^{-1/2}\sup_{s\leq t}\|\Upsilon^n(s)\| \to 0$ in distribution for any $t$, and therefore $n^{-1/2}\Upsilon^n \Rightarrow 0$. As a result, $n^{-1/2}\hat{X}^n \Rightarrow 0$, which implies that $\bar{X}^n \Rightarrow \rho$. Using $\mathbb{1}\cdot\rho = 1$ and $\mathbb{1}\cdot\bar{\Phi}^n = (\mathbb{1}\cdot\bar{X}^n - 1)^+$, we have that $\mathbb{1}\cdot\bar{\Phi}^n \Rightarrow 0$. Now $\bar{\Phi}_i^n \Rightarrow 0$ follows since $\bar{\Phi}_i^n \geq 0$. Using $\bar{X}^n = \bar{\Phi}^n + \bar{\Psi}^n$, we have that $\bar{\Psi}^n \Rightarrow \rho$.

(iii) By (i), $\hat{A}^n \Rightarrow A$. By (i) and (ii) and a time change lemma (cf. [4]), it directly follows that $\hat{S}_i^n(\int_0^\cdot \bar{\Psi}_i^n(s)\,ds) \Rightarrow S_i(\rho_i\cdot)$. A use of (i), (ii) and a time change lemma also shows that $\hat{R}_i^n(\int_0^\cdot \bar{\Phi}_i^n(s)\,ds) \Rightarrow 0$. Hence by (14),

(56) $$\hat{W}_i^n \Rightarrow r^{-1}rW = W,$$

where $W$ is a standard $k$-dimensional Brownian motion.

Since $\hat{W}^n$ are relatively compact, they are tight. Hence by [4], Theorem 16.8, for each $t$, $\lim_{m\to\infty}\limsup_{n\to\infty} P(\|\hat{W}^n\|_t \geq m) = 0$. By (23) and the Lipschitz property of the functions $x \mapsto b^n(x, u)$, uniformly in $x$, $u$ and $n$,

$$\|\hat{X}^n(t)\| \leq \|\hat{X}^{0,n}\| + \|\hat{W}^n(t)\| + c\int_0^t (1 + \|\hat{X}^n(s)\|)\,ds.$$

By Gronwall's inequality, using the boundedness of $\hat{X}^{n,0}$, $n \in \mathbb{N}$, we have

(57) $$\|\hat{X}^n\|_t \leq ce^{ct}(1 + \|\hat{W}^n\|_t).$$

This shows that, for each $t$,

(58) $$\lim_{m\to\infty}\limsup_{n\to\infty} P(\|\hat{X}^n\|_t \geq m) = 0.$$

Fix $T$. It follows from (23) that, for any $s,t \in [0,T]$ with $s < t$,

(59) $$\|\hat{X}^n(t) - \hat{X}^n(s)\| \leq \|\hat{W}^n(t) - \hat{W}^n(s)\| + c\int_s^t (1 + \|\hat{X}^n(s)\|)\,ds$$
$$\leq \|\hat{W}^n(t) - \hat{W}^n(s)\| + c(t-s)(1 + \|\hat{X}^n\|_T).$$

Recall the modulus of continuity defined for $x \in \mathbb{D}(\mathbb{R}^k)$ restricted to $[0,T]$ (cf. [4], page 171) as

$$w'_T(x,\delta) = \inf \max_{1 \leq i \leq v} w(x, [t_{i-1}, t_i)),$$

where the infimum is taken over all decompositions $[t_{i-1}, t_i)$, $1 \leq i \leq v$, of $[0,T)$ such that $t_i - t_{i-1} > \delta$ for $1 \leq i \leq v$. Here, for $S \subset [0,T)$,

$$w(x, S) = \sup_{s,t \in S} \|x(s) - x(t)\|.$$

By tightness of $\hat{W}^n$, Theorem 16.8 of [4] implies that, for each $t$ and $\varepsilon$,

$$\lim_{\delta \to 0}\limsup_{n\to\infty} P(w'_t(\hat{W}^n, \delta) \geq \varepsilon) = 0.$$

Using (59), a similar statement follows for $\hat{X}^n$, namely that, for each $t \leq T$ and $\varepsilon$,

(60) $$\lim_{\delta \to 0}\limsup_{n\to\infty} P(w'_t(\hat{X}^n, \delta) \geq \varepsilon) = 0.$$

By (59) and (60), and since $T$ is arbitrary, the tightness of $\hat{X}^n$ follows from Theorem 16.8 of [4].

Noting that $\|Y^n(t)\| \leq ct(1 + \|\hat{X}^n\|_t)$, and $\|Z^n(t)\| \leq ct(1 + \|\hat{X}^n\|_t^m)$ ($m$ of the $L$), and that for $s,t \leq T$, $\|Y^n(t) - Y^n(s)\| \leq c|t-s|(1 + \|\hat{X}^n\|_T)$ and $\|Z^n(t) - Z^n(s)\| \leq c|t-s|(1 + \|\hat{X}^n\|_T^m)$, the tightness of $Y^n$ and of $Z^n$ follows from (58) using again Theorem 16.8 of [4]. □

We use the following (very special case of a) result of Kurtz and Protter [28]. Let $(F_t)$ be a filtration. A cadlag, $(F_t)$-adapted process $V$ is a semimartingale if $V = M + N$, where $M$ is an $(F_t)$-local martingale, and the paths of $N$ are of bounded variation over finite time intervals. An $\mathbb{R}^k$-valued process is an $(F_t)$-semimartingale if each component is a semimartingale. Write $\int U\, dV$ for $\int_0^\cdot U(s-) \cdot dV(s)$. A cadlag process $V$ has bounded jumps if there is a constant $c$ such that $\|V(s) - V(s-)\| \leq c$, $s \in (0, \infty)$, a.s. Denote by $[M]$ the quadratic variation process associated with $M$, and by $T_t(N)$ the total variation of $N$ over $[0, t]$.

LEMMA 5. *For each $n$, let $(U^n, V^n)$ be an $(F_t^n)$-adapted process with sample paths in $\mathbb{D}((\mathbb{R}^k)^2)$ and let $V^n$ be an $(F_t^n)$-semimartingale with bounded jumps. Let $V^n = M^n + N^n$ be a decomposition of $V^n$ into an $(F_t^n)$-local martingale and a process with finite variation. Suppose*

(61) $$\text{for each } t > 0, \qquad \sup_n E[[M^n]_t + T_t(N^n)] < \infty.$$

*If $(U^n, V^n) \Rightarrow (U, V)$ in the Skorohod topology on $\mathbb{D}((\mathbb{R}^k)^2)$, then $V$ is a semimartingale with respect to a filtration to which $U$ and $V$ are adapted, and $(U^n, V^n, \int U^n\, dV^n) \Rightarrow (U, V, \int U\, dV)$ in the Skorohod topology on $\mathbb{D}((\mathbb{R}^k)^3)$.*

PROOF. The proof follows from Theorem 2.2 of [28] on taking, for $\alpha > 0$, $\tau_n^\alpha = \alpha + 1$, noting that $V_n^\delta = V_n$ if $\delta$ is a fixed large constant. $\square$

LEMMA 6. *Let the assumptions of Theorem 4(i) hold. Denote by $(X, Y, Z, W)$ a limit point of $(\hat{X}^n, Y^n, Z^n, \hat{W}^n)$ along a subsequence. Let $(F_t)$ denote the filtration generated by $(X, Y, W)$. Then $W$ is an $(F_t)$-standard Brownian motion, $X$, $Y$ and $Z$ have continuous sample paths, and $Y$ has sample paths of bounded variation over finite time intervals. Moreover, $\int e^{-\gamma s} Df(\hat{X}_s^n) \cdot dY_s^n \Rightarrow \int e^{-\gamma s} Df(X_s) \cdot dY_s$ along the subsequence, where $f$ is the solution to (41).*

PROOF. The processes $Y$ and $Z$ have continuous sample paths since $Y^n$ and $Z^n$ do (see Theorem 3.10.2(a) of [10]). Since $\hat{X}^n = \hat{X}^{0,n} + r\hat{W}^n + Y^n$, and $\hat{W}^n$ converges in distribution to a Brownian motion [cf. (56)], $X = x + rW + Y$ has continuous sample paths. To see that $Y$ has sample paths of bounded variation, write $Y^n = Y^{n,+} - Y^{n,-}$, where $Y_i^{n,+}(t) = \int_0^t (\dot{Y}_i^n(s))^+\, ds$, $Y_i^{n,-}(t) = \int_0^t (\dot{Y}_i^n(s))^-\, ds$. By definition (50) of $Y_t^n$ and (24) of $b^n$,

(62) $$Y^{n,+}(t) \vee Y^{n,-}(t) \leq c \int_0^t (1 + \|\hat{X}_s^n\|)\, ds,$$
$$(Y^{n,+}(t) - Y^{n,+}(s)) \vee (Y^{n,-}(t) - Y^{n,-}(s)) \leq c|t - s|(1 + \|\hat{X}^n\|_t),$$

where $c$ does not depend on $t, n$. Thus it follows from the tightness of $\hat{X}^n$ that $(Y^{n,+}, Y^{n,-})$ is tight. Let $(Y^+, Y^-)$ denote any subsequential limit point in $(\mathbb{D}(\mathbb{R}^k))^2$. Since $Y^{n,+}$ and $Y^{n,-}$ have continuous sample paths, so do $Y^+$ and $Y^-$, and therefore $Y = Y^+ - Y^-$. Since $Y^+$ and $Y^-$ have nondecreasing sample paths, $Y$ has sample paths of bounded variation over $[0, t]$ for any $t$.

Next we apply Lemma 5 with $U^n = e^{-\gamma t} Df(\hat{X}^n(t))$, $V^n = Y^n$, and $(F_t^n) = (\mathcal{F}_t^n)$ of (9). By Definition 2 and the definition of $Y^n$, clearly $\hat{X}^n$ and $Y^n$ are adapted to $(F_t^n)$. We decompose $Y^n = M^n + N^n$ as $M^n = 0$, $N^n = Y^n$. By (62), and Lemma 3, (61) holds. By the continuous mapping theorem, $(e^{-\gamma t} Df(\hat{X}^n(t)), Y^n(t))$ converges to $(e^{-\gamma t} Df(X(t)), Y(t))$ in the Skorohod topology on $(\mathbb{D}(\mathbb{R}^k))^2$. By continuity of the sample paths of $Y^n$, it follows that the convergence in fact holds in the Skorohod topology on $\mathbb{D}((\mathbb{R}^k)^2)$ (see Proposition 6.3.2 of [10]). As a result of Lemma 5, $\int e^{-\gamma t} Df(\hat{X}^n(t)) \cdot dY^n(t) \Rightarrow \int e^{-\gamma t} Df(X(t)) \cdot dY(t)$.

It was shown in the proof of Lemma 4 [cf. (56)] that $\hat{W}^n$ converges to a standard Brownian motion. To see that $W$ is in fact an $(F_t)$-Brownian motion, note that by definition it is adapted to $(F_t)$. It remains to show that, for each $t$, $F_t$ is independent of $\sigma\{W_{t+u} - W_t : u > 0\}$. Fix $t \geq 0$, $u \geq 0$ and $0 \leq s \leq t$. Write $\alpha^n = (\hat{X}_s^n, Y_s^n, \hat{W}_s^n)$ and $\alpha = (X_s, Y_s, W_s)$. By (14), using the notation (8), and denoting

$$\tilde{S}_i^n = \hat{S}_i^n(n^{-1} T_i^n(t)) - \hat{S}_i^n(n^{-1} T_i^n(t+u)),$$
$$\tilde{R}_i^n = \hat{R}_i^n(n^{-1} \overset{\circ}{T}_i^n(t)) - \hat{R}_i^n(n^{-1} \overset{\circ}{T}_i^n(t+u)),$$

we have

$$r_i(\hat{W}_i^n(t+u) - \hat{W}_i^n(t)) = \hat{A}_i^n(t+u) - \hat{A}_i^n(t) - \tilde{S}_i^n - \tilde{R}_i^n$$
$$= r_i \beta_i^n + \delta_i^n,$$

where

$$r_i \beta_i^n = \hat{A}_i^n(\tau_i^n(t) + u) - \hat{A}_i^n(\tau_i^n(t)) - \tilde{S}_i^n - \tilde{R}_i^n$$

and

$$\delta_i^n = \hat{A}_i^n(t+u) - \hat{A}_i^n(t) - \hat{A}_i^n(\tau_i^n(t) + u) + \hat{A}_i^n(\tau_i^n(t)).$$

Let $f : (\mathbb{R}^k)^3 \to \mathbb{R}$ and $g : \mathbb{R}^k \to \mathbb{R}$ be bounded continuous. By (9) and (10), $\beta^n$ is measurable on $\mathcal{G}_t^n$ and $\alpha^n$ is measurable on $\mathcal{F}_t^n$. By the admissibility assumption and Definition 2,

(63) $$Ef(\alpha^n) g(\beta^n) = Ef(\alpha^n) Eg(\beta^n).$$

Since $\tau_i^n(t)$ converges in distribution to zero, and $\hat{A}_i^n$ converges in distribution to a continuous process, it follows by a random change of time lemma ([4], page 151) that $\delta_i^n$ converges in distribution to zero. As a result, $\beta^n$ converges in distribution to $W_{t+u} - W_t$. Using (63), the convergence $(\hat{X}^n, Y^n, \hat{W}^n) \Rightarrow (X, Y, W)$ and the continuous mapping theorem, it follows that

$$(64) \qquad Ef(\alpha)g(W_{t+u} - W_t) = Ef(\alpha)Eg(W_{t+u} - W_t).$$

By approximating indicator functions of closed sets of $(\mathbb{R}^k)^3$ (and resp. $\mathbb{R}^k$) by continuous functions $f$ (resp. $g$), it follows that (64) holds when $f$ and $g$ are replaced by such indicator functions. Since $u \geq 0$ and $s \leq t$ are arbitrary, an application of the Dynkin class theorem (Theorem 1.4.2 of [9]) shows that $F_t$ is independent of $\sigma\{W_{t+u} - W_t : u > 0\}$. Since also $t$ is arbitrary, it follows that $W$ is an $(F_t)$-Brownian motion. $\square$

PROOF OF THEOREM 4. We first prove part (ii). Recall that (49) holds. Let $(X, Y, Z, W)$ be a weak limit point of $(\hat{X}^n, Y^n, Z^n, \hat{W}^n)$ and let $(F_t)$ be the filtration generated by $(X, Y, W)$. By Lemma 6, $X_t = x + rW_t + Y_t$, $W$ is a standard $(F_t)$-Brownian motion and the sample paths of $Y$ have bounded variation over finite time intervals. Just as before, an application of Itô's formula and the fact that $f$ satisfies the HJB equation (41) give

$$(65) \quad \begin{aligned} e^{-\gamma t} f(X_t) &= f(x) + \int_0^t e^{-\gamma s} Df(X_s) \cdot r \, dW_s \\ &\quad + \int_0^t e^{-\gamma s} Df(X_s) \cdot dY_s - \int_0^t e^{-\gamma s} H(X_s, Df(X_s)) \, ds. \end{aligned}$$

By (48),

$$(66) \quad \begin{aligned} \int_0^t e^{-\gamma s} K_s^n \, ds &= \int_0^t e^{-\gamma s} Df(\hat{X}_s^n) \cdot dY_s^n \\ &\quad + e_n(t) + Z_t^n - \int_0^t e^{-\gamma s} H(\hat{X}_s^n, Df(\hat{X}_s^n)) \, ds, \end{aligned}$$

where

$$e_n(t) = \int_0^t e^{-\gamma s}(b(\hat{X}_s^n, u_s^n) - b^n(\hat{X}_s^n, u_s^n)) \cdot Df(\hat{X}_s^n) \, ds.$$

By definition of the functions $b$ and $b^n$ and by Assumption 1,

$$\|b(\hat{X}_s^n, u_s^n) - b^n(\hat{X}_s^n, u_s^n)\| \leq \varepsilon_n(1 + \|\hat{X}_s^n\|),$$

where $\varepsilon_n \to 0$. Therefore, Lemma 3 and the continuous mapping theorem imply that $e_n \Rightarrow 0$. We get from (49) and Lemma 6, using continuity of $x \mapsto H(x, Df(x))$,

$$(67) \qquad \int_0^t e^{-\gamma s} Df(X_s) \cdot dY_s + Z_t - \int_0^t e^{-\gamma s} H(X_s, Df(X_s)) \, ds = 0.$$

Combining (65) and (67),

$$0 \le e^{-\gamma t} f(X_t) = f(x) + \int_0^t e^{-\gamma s} Df(X_s) \cdot r \, dW_s - Z_t.$$

Hence

(68) $$\forall t, \qquad EZ_t \le f(x).$$

Fix an arbitrary $\delta > 0$. By Lemma 3 and Assumption 2, there is $T$ such that

(69) $$E \int_T^\infty e^{-\gamma s} L(\hat{X}_s^n, u_s^n) \, ds \le \delta$$

for all $n$. Since $Z^n \Rightarrow Z$ and $Z$ has continuous sample paths, $Z_T^n$ converges in distribution to $Z_T$. By Jensen's inequality, Assumption 2 and Lemma 3,

(70) $$E(Z_T^n)^{1+\varepsilon/m_L} \le cE \int_0^T e^{-\gamma(1+\varepsilon/m_L)s}(1 + \|\hat{X}_s^n\|^{m_L+\varepsilon}) \, ds$$
$$\le c,$$

where $c$ does not depend on $n$. Hence $Z_T^n$, $n \in \mathbb{N}$, are uniformly integrable, and one has $EZ_T^n \to EZ_T$ as $n \to \infty$. By (68) and (69), we therefore have that

$$\limsup_{n \to \infty} E \int_0^\infty e^{-\gamma s} L(\hat{X}_s^n, u_s^n) \, ds \le f(x) + \delta.$$

Since $\delta > 0$ is arbitrary, it can be dropped, and part (ii) of the theorem follows.

Next we prove part (i). Arguing as in part (ii) but using $K_t^n \ge 0$ instead of (49), we have that (65) holds and

$$\int_0^t e^{-\gamma s} Df(X_s) \cdot dY_s + Z_t - \int_0^t e^{-\gamma s} H(X_s, Df(X_s)) \, ds \ge 0.$$

Hence

(71) $$e^{-\gamma t} f(X_t) \ge f(x) + \int_0^t e^{-\gamma s} Df(X_s) \cdot r \, dW_s - Z_t.$$

By Proposition 5 and Lemma 3,

$$Ef(\hat{X}_t^n) \le c(1 + t^{m_L}),$$

for $t \ge 0$ and $n \in \mathbb{N}$. Since for each $t$, $f(\hat{X}_t^n)$ converges in distribution to $f(X_t)$, and $f(\hat{X}_t^n)$ are uniformly integrable [arguing as in (70), using the growth condition of Proposition 5(i)], one has $Ef(X_t) \le c(1 + t^{m_L})$, where again $c$ does not depend on $t$. We therefore have, from (71),

$$EZ_t \ge f(x) - \alpha(t),$$

where $\alpha(t) \to 0$ as $t \to \infty$. Note that as in part (ii), given $\delta > 0$, (69) holds for all $T$ large enough, and that, for each $T$, $Z_T^n$, $n \in \mathbb{N}$, are uniformly integrable. Hence

$$\liminf_{n \to \infty} E \int_0^\infty L(\hat{X}_s^n, u_s^n) \, ds \geq E Z_T - \delta \geq f(x) - \alpha(T) - \delta.$$

Part (i) of the theorem now follows on taking $T \to \infty$ and $\delta \to 0$. $\square$

PROOF OF THEOREM 2. We only need to show that the proposed SCPs satisfy the conditions of Theorem 4(ii). The work-conservation condition (19) holds for both of the proposed SCPs, by definition. To conclude parts (i) and (ii), it remains to show that in both cases (49) holds. Part (iii) is treated thereafter.

(i) *The P-SCP.* Fix $T$. Let $\Omega^n$ denote the event that (31) is met for all $t \in [0, T]$. Recall that on $\Omega^n$, the P-SCP sets

$$\Phi(t) = \Theta[(\mathbb{1} \cdot X^n(t) - n)^+ h(\hat{X}^n(t))], \qquad t \in [0, T].$$

Let

(72) $$U^n = (\mathbb{1} \cdot \hat{X}^n)^+ h(\hat{X}^n), \qquad V^n = \hat{\Phi}^n - U^n.$$

Recall that $h$ satisfies

$$b(x, h(x)) \cdot Df(x) + L(x, h(x)) = H(x, Df(x)), \qquad x \in \mathbb{R}^k.$$

Note that for $x$ with $\mathbb{1} \cdot x \leq 0$, $b(x, u)$ is independent of $u$ [see (26)] and so is $L(x, u) = \tilde{L}((\mathbb{1} \cdot x)^+ u, x - (\mathbb{1} \cdot x)^+ u)$ [see (21)]. Hence

(73) $$\inf_{u \in \mathbb{S}^k} [b(x, u) \cdot p + L(x, u)] = b(x, v) \cdot p + L(x, v), \qquad v \in \mathbb{S}^k, \quad \mathbb{1} \cdot x \leq 0.$$

For $t$ such that $\mathbb{1} \cdot \hat{X}_t^n \leq 0$, (73) and (48) imply that $K_t^n = 0$. Next consider $t$ such that $\mathbb{1} \cdot \hat{X}_t^n > 0$. We have

$$u_t^n = \hat{\Phi}^n(t)(\mathbb{1} \cdot \hat{X}_t^n)^{-1} = h(\hat{X}_t^n) + V_t^n (\mathbb{1} \cdot \hat{X}_t^n)^{-1}.$$

By assumption, $\tilde{L}$ is uniformly continuous on compacts. For each $\kappa$, let $\alpha^\kappa(\delta)$ be such that $\alpha^\kappa(\delta) \downarrow 0$ as $\delta \downarrow 0$, and $|\tilde{L}(\phi, \psi) - \tilde{L}(\phi', \psi')| \leq \alpha^\kappa(\delta)$ whenever $\|\phi\|, \|\phi'\|, \|\psi\|, \|\psi'\| \leq \kappa$, and $\|\phi - \phi'\| \vee \|\psi - \psi'\| \leq \delta$. Then using (21), writing $\xi_t^n = \mathbb{1} \cdot \hat{X}_t^n$, the following holds on the event $\Omega^{n,\kappa} := \Omega^n \cap \{\|\hat{\Phi}^n\|_T^* + \|\hat{\Psi}^n\|_T^* + \|\hat{X}^n\|_T^* \leq \kappa\}$:

(74)
$$\begin{aligned}|K_t^n| &= |(b(\hat{X}_t^n, u_t^n) - b(\hat{X}_t^n, h(\hat{X}_t^n))) \cdot Df(\hat{X}_t^n) \\ &\quad + \tilde{L}((\xi_t^n)^+ u_t^n, \hat{X}_t^n - (\xi_t^n)^+ u_t^n) \\ &\quad - \tilde{L}((\xi_t^n)^+ h(\hat{X}_t^n), \hat{X}_t^n - (\xi_t^n)^+ h(\hat{X}_t^n))| \\ &\leq c\|V_t^n\| \|Df(\hat{X}_t^n)\| + \alpha^\kappa(c\|V_t^n\|).\end{aligned}$$

By (30), $\|V_t^n\| \le 2kn^{-1/2}$. As a result, $|K^n|_T^* \le \varepsilon_n$ on $\Omega^{n,\kappa}$ for some $\varepsilon_n \to 0$. The events $\Omega^n$ have probability tending to 1 as $n \to \infty$, as follows from the convergence $\bar{X}^n \Rightarrow \rho$ shown in Lemma 4. The tightness of $\hat{X}^n$ (see Lemma 4), (19) and the fact that $\hat{\Phi}_t^n \in \mathbb{R}_+^k$ imply that

(75) $$\lim_{\kappa \to \infty} \liminf_{n \to \infty} P(\Omega^{n,\kappa}) = 1.$$

Therefore $|K^n|_T^*$ converges to zero in distribution. Since $T$ is arbitrary, $K^n \Rightarrow 0$, and (49) holds.

(ii) *The N-SCP* (i). Fix $T$. Let $U^n$ and $V^n$ be defined as in (72). A review of the previous paragraph shows that, replacing throughout $\Omega^n$ by $\Omega$, (74) and (75) still hold. Fix $\varepsilon_0 > 0$. We next estimate, for any $\varepsilon > 0$,

$$\limsup_n P\Big(\sup_{t \in [\varepsilon_0, T]} \|V^n(t)\| > 8k\varepsilon\Big).$$

Fix $i \in K$. If either $V_i^n(t) < 0$ or $\hat{\Phi}_i^n(t) = 0$ holds for all $t \in [s, r)$, then within this time interval, the SCP does not route any class-$i$ customer to the service pool. Therefore by (4), for $t \in [s, r)$,

(76) $$\hat{\Phi}_i^n(t) = \hat{\Phi}_i^n(s) + n^{-1/2} A_i^n(s, t) - n^{-1/2} \Delta_i^n(s, t),$$

where we write

$$A_i^n(s, t) = A_i^n(t) - A_i^n(s),$$

$$\Delta_i^n(s, t) = R_i^n\Big(\int_0^t \Phi_i^n(z)\,dz\Big) - R_i^n\Big(\int_0^s \Phi_i^n(z)\,dz\Big).$$

Given $\varepsilon > 0$,

(77) $$P\Big(\inf_{t \in [\varepsilon_0, T]} V_i^n(t) < -4\varepsilon\Big) \le P((\Omega^{n,\kappa})^c) + P(\Omega_1^{n,k}) + P(\Omega_2^{n,\kappa}),$$

where

$$\Omega_1^{n,\kappa} = \Omega^{n,\kappa} \cap \Big\{ \exists\, \varepsilon_0 \le s \le r \le T : V_i^n(s) \ge -\varepsilon,$$

$$\sup_{t \in [s, r)} V_i^n(t) \le -\varepsilon,\ V_i^n(r) \le -4\varepsilon \Big\},$$

$$\Omega_2^{n,\kappa} = \Omega^{n,\kappa} \cap \Big\{ \sup_{t \in [0, \varepsilon_0]} V_i^n(t) < -\varepsilon \Big\}.$$

Using the local Hölder property of $h$ on $\mathbb{X}$, for any $\kappa$, there are $c_\kappa > 0$ and $p_\kappa \in (0, 1]$ such that, on $\Omega^{n,\kappa}$,

$$|U_i^n(t) - U_i^n(s)| \le c_\kappa \|\hat{X}^n(t) - \hat{X}^n(s)\|^{p_\kappa}$$
$$+ (\varepsilon/4)\mathbf{1}_{\{\mathbb{1}\cdot \hat{X}^n(s) < \varepsilon/8\}} + (\varepsilon/4)\mathbf{1}_{\{\mathbb{1}\cdot \hat{X}^n(t) < \varepsilon/8\}}$$
$$\le c_\kappa \|\hat{X}^n(t) - \hat{X}^n(s)\|^{p_\kappa} + \varepsilon/2.$$

Writing
$$n^{-1/2}\Delta_i^n(s,t)$$
$$= \hat{R}_i^n\left(\int_0^t \bar{\Phi}_i^n(z)\,dz\right) - \hat{R}_i^n\left(\int_0^s \bar{\Phi}_i^n(z)\,dz\right) + n^{1/2}\theta_i \int_s^t \bar{\Phi}_i^n(z)\,dz,$$
and using $\|\bar{\Phi}^n\|_T^* \leq \kappa n^{-1/2}$ on $\Omega^{n,\kappa}$, we have
$$n^{-1/2}|\Delta_i^n(s,t)| \leq 2\|\hat{R}^n\|_{\kappa T n^{-1/2}}^* + c\kappa(t-s).$$
Hence on $\Omega_1^{n,\kappa}$, for $t \in [s,r]$,
$$V_i^n(t) = V_i^n(s) + (\hat{\Phi}_i^n(t) - \hat{\Phi}_i^n(s)) - (U_i^n(t) - U_i^n(s))$$
$$\geq -\varepsilon - \varepsilon/2 + n^{-1/2}A_i^n(s,t)$$
$$- c_\kappa \|\hat{X}^n(t) - \hat{X}^n(s)\|^{p_\kappa} - 2\|\hat{R}^n\|_{\kappa T n^{-1/2}}^* - c\kappa(t-s).$$
On $\Omega_1^{n,\kappa}$ we also have $V_i^n(r) \leq -4\varepsilon$. Let $\beta > 0$ and write $\tilde{\beta} = \beta + c\kappa$. Therefore
$$(78) \qquad P(\Omega_1^{n,\kappa}) \leq P(\Omega_{1,1}^{n,\kappa}) + P(\Omega_{1,2}^{n,\kappa}) + P(\Omega_{1,3}^{n,\kappa}),$$
where
$$\Omega_{1,1}^{n,\kappa} = \{\exists\, 0 \leq s \leq r \leq T : n^{-1/2}A^n(s,r) \leq -\varepsilon + \tilde{\beta}(r-s)\},$$
$$\Omega_{1,2}^{n,\kappa} = \{\exists\, 0 \leq s \leq r \leq T : c_\kappa \|\hat{X}^n(r) - \hat{X}^n(s)\|^{p_\kappa} \geq \varepsilon + \beta(r-s)\},$$
$$\Omega_{1,3}^{n,\kappa} = \{2\|\hat{R}^n\|_{\kappa T n^{-1/2}}^* \geq \varepsilon/2\}.$$
Using the monotonicity of the processes $A_i^n$ and the uniform convergence of $n^{-1}A_i^n$ on $[0,T]$ to $\tilde{A}_i(t) = \lambda_i t$, as follows from the convergence $\hat{A}^n \Rightarrow A$ [see Lemma 4(i)],
$$\limsup_n P(\Omega_{1,1}^{n,\kappa})$$
$$(79) \qquad \leq \limsup_n P(\exists\, 0 \leq s \leq r \leq T : r - s \geq \varepsilon/\tilde{\beta},\ n^{-1}A^n(s,r) \leq n^{-1/2}\tilde{\beta}T)$$
$$\leq \limsup_n P\left(\sup_{t \leq T} \|n^{-1}A^n(t) - \lambda_i t\| \geq c\right)$$
$$= 0.$$
Also,
$$P(\Omega_{1,2}^{n,\kappa}) \leq P(\exists\, 0 \leq s \leq r \leq T, r-s > \beta^{-1/2} : c_\kappa \|\hat{X}^n(r) - \hat{X}^n(s)\|^{p_\kappa} \geq \beta^{1/2})$$
$$+ P(\exists\, 0 \leq s \leq r \leq T, r-s \leq \beta^{-1/2} : c_\kappa \|\hat{X}^n(r) - \hat{X}^n(s)\|^{p_\kappa} \geq \varepsilon)$$
$$\leq P(2c_\kappa^{p_\kappa^{-1}}\|\hat{X}^n\|_T^* \geq \beta^{1/(2p_\kappa)}) + P(w(\hat{X}^n|_{[0,T]}, \beta^{-1/2}) \geq (\varepsilon/c_\kappa)^{p_\kappa^{-1}}).$$

By Lemma 6, the processes $\hat{X}^n$ are tight and converge along subsequences to processes with continuous sample paths. Therefore

$$\lim_{\beta \to \infty} \limsup_n P(\Omega_{1,2}^{n,\kappa}) = 0. \tag{80}$$

The convergence of $\hat{R}^n$ to a Brownian motion (Lemma 4) implies

$$\lim_n P(\Omega_{1,3}^{n,\kappa}) = 0. \tag{81}$$

By a similar argument, on $\Omega_2^{n,\kappa}$,

$$V_i^n(0) + n^{-1/2} A_i^n(\varepsilon_0) - c_\kappa \|\hat{X}^n(\varepsilon_0) - \hat{X}^n(0)\|^{p_\kappa} - \varepsilon/2 - 2\|\hat{R}^n\|^*_{\kappa T n^{-1/2}} - c\kappa \varepsilon_0$$
$$\leq V_i^n(\varepsilon_0) \leq -\varepsilon.$$

Hence, for some constant $c'_\kappa$,

$$\lim_n P(\Omega_2^{n,\kappa}) \leq \lim_n P(n^{-1} A_i^n(\varepsilon_0) \leq c'_\kappa n^{-1/2} + 2n^{-1/2}\|\hat{R}^n\|^*_{\kappa T n^{-1/2}}) \tag{82}$$
$$= 0,$$

where the last equality follows from the convergence in distribution of $n^{-1} A_i^n(\varepsilon_0)$ to $\lambda_i \varepsilon_0$ and of $\hat{R}^n$ to a Brownian motion. Combining (77)–(82) shows that

$$\limsup_n P\left(\inf_{t \in [\varepsilon_0, T]} V_i^n(t) < -4\varepsilon\right) \leq \limsup_n P((\Omega^{n,\kappa})^c).$$

Note that by (19), $\mathbb{1} \cdot V^n = 0$. Hence $\|V^n\| = 2(\mathbb{1} \cdot V^n)^-$. Since $i \in K$ is arbitrary, it follows that

$$\limsup_n P\left(\sup_{t \in [\varepsilon_0, T]} \|V^n(t)\| > 8k\varepsilon\right) \leq \limsup_n P((\Omega^{n,\kappa})^c).$$

Combining this with (74) (assuming without loss that, for each $\kappa$, $\alpha^\kappa$ is bounded), (75) and the fact that $\varepsilon, \varepsilon_0 > 0$ and $T$ are arbitrary, it follows that $\int_0^{\cdot} e^{-\gamma s} K_s^n \, ds \Rightarrow 0$. Therefore (49) holds and this concludes the proof that both SCPs satisfy the conditions of Theorem 4(i). Parts (i) and (ii) of Theorem 2 follows.

(iii) *The N-SCP* (ii). To prove this part, it suffices to show that, for each $\delta$, there is a locally Lipschitz $h'$ such that the N-SCP (i) applied to $h'$ gives

$$\limsup_{n \to \infty} E \int_0^\infty e^{-\gamma t} L(\hat{X}_t^n, u_t^n) \, dt \leq V(x) + \delta. \tag{83}$$

Recall from the proof of Theorem 3 that $h : \mathbb{R}^k \to \mathbb{S}^k$ is a function satisfying

$$\varphi(x, h(x)) = \inf_{v \in \mathbb{S}^k} \varphi(x, v) =: \varphi^*(x),$$

where

$$\varphi(x,u) = b(x,u) \cdot Df(x) + L(x,u).$$

For each $\varepsilon > 0$, let $h^\varepsilon : \mathbb{R}^k \to \mathbb{S}^k$ be a function defined as

$$h^\varepsilon(x) = \frac{\sum_y d(x,y,\varepsilon) h(y)}{\sum_y d(x,y,\varepsilon)},$$

where both sums extend over $y \in \varepsilon \mathbb{Z}^k \cap B(x, \varepsilon k^{1/2})$, and $d(x,y,\varepsilon)$ denotes the Euclidean distance from $y$ to the boundary $\partial B(x, \varepsilon k^{1/2})$. It is easy to check that $h^\varepsilon$ is locally Lipschitz. Write $\tilde{d}(x,y,\varepsilon) = d(x,y,\varepsilon)/\sum_{y'} d(x,y',\varepsilon)$. By assumption, $u \mapsto L(x,u)$ is convex, and since $u \mapsto b(x,u)$ is affine, $u \mapsto \varphi(x,u)$ is convex. Using Jensen's inequality, uniform continuity of $(x,u) \mapsto \varphi(x,u)$ and of $x \mapsto \varphi^*(x)$ on compacts, for each $\delta > 0$, there is $\varepsilon$ such that

$$\begin{aligned}
\varphi(x, h^\varepsilon(x)) &= \varphi\left(x, \sum_y \tilde{d}(x,y,\varepsilon) h(y)\right) \\
&\leq \sum_y \tilde{d}(x,y,\varepsilon) \varphi(x, h(y)) \\
&\leq \sum_y \tilde{d}(x,y,\varepsilon) \varphi(y, h(y)) + \delta/2 \\
&= \sum_y \tilde{d}(x,y,\varepsilon) \varphi^*(y) + \delta/2 \\
&\leq \varphi^*(x) + \delta.
\end{aligned}$$

Everywhere in the above display, the sum extends over $y \in \varepsilon \mathbb{Z}^k \cap B(x, \varepsilon k^{1/2})$. A review of the proof of Theorems 2 and 4 shows that, upon applying N-SCP (i) with $h^\varepsilon$, (83) holds. By taking an appropriate sequence $h_n = h^{\varepsilon_n}$, it is then clear that N-SCP (ii) applied to $h_n$ admits the conclusion of Theorem 4, and therefore (32). □

## 5. Further research.

5.1. *Work-encouraging SCPs.* We have restricted our analysis to work-conserving SCPs. However, our results regarding asymptotic optimality among *all* admissible SCPs hold, in fact, under the additional condition that the cost functions are work encouraging (cf. Definition 6). Recall that with each admissible SCP we have associated a cost of the form [cf. (20)]

$$C^n = E \int_0^\infty e^{-\gamma t} \tilde{L}(\hat{\Phi}^n(t), \hat{\Psi}^n(t)) \, dt.$$

DEFINITION 6. We say that the cost function $\tilde{L}$ (or the corresponding cost function $L$) is *work encouraging* if, for each $n$, the infimum of $C^n$ over all admissible SCPs is equal to that over all work-conserving admissible SCPs.

COROLLARY 1. *Let all assumptions of Theorem 2 apply, except the assumption that the SCPs are work conserving. Then the conclusions of Theorem 2 prevail, given that the cost function $\tilde{L}$ is work encouraging.*

Although in many cases it is intuitively clear that work conservation is optimal (for P-SCPs, not for N-SCPs), in the presence of abandonments, and in the generality of our setting, this turns out to be nontrivial to prove. We intend to treat the issue in a future work. We end this section with a few examples that are intended to exhibit some of the subtleties of this point, and to indicate how it can be dealt with. The arguments should be considered as proof outlines only.

First, consider the expected discounted number of customers of a particular class, say class 1, present in the system. If $\theta_1 \leq \mu_1$, then class-1 customers leave the system faster when they are served than when they are in the queue. Hence a good policy will attempt to serve these customers as much as possible, and will be work conserving. On the other hand, if $\theta_i > \mu_i$, then customers leave the system by abandoning the queue faster than by being served, and as a result, a policy which minimizes the cost will not schedule any services at all.

More subtle are the costs associated with queue length and abandonment. We argue heuristically that if $\theta_1 > \mu_1$, then there are cost functions $\ell$, nondecreasing as a function of $\hat{\Phi}_i$ for each $i$, for which work conservation is actually not optimal. Suppose that $\theta_1 > \mu_1$, and the cost is 0 for $\hat{\Phi}_1 \leq c$, and 1 for $\hat{\Phi}_1 > c$, where $c > 0$. If $\hat{\Phi}_1 \leq c$, no cost is incurred, and customers leave faster if in the queue than if in service. Thus an SCP that keeps customers in the queue would do better than a work-conserving SCP.

Consider the case where $\ell$ is linear in $\hat{\Phi}_i$, and $\theta_i, \mu_i$ are arbitrary. We argue that work-conserving policies are optimal. We use coupling. A sample path is considered under an SCP that is not always work conserving. The coupling is used to show that if the SCP is changed to be work conserving, the cost will be no higher than for the original SCP. In view of the discussion on costs of abandonment, one can use the relation between abandonment rate and expected queue length to obtain the result. Consider a sample path under an SCP that leaves customers in the queue when there are idle servers. Modify it by moving a customer into service. Keep that customer in service until the earliest of: (i) it completes its service, (ii) its "twin" (i.e., the customer in the original system that is in the queue) abandons, or (iii) the original SCP needs to use the server. In cases (i) and (ii), the cost of the modified SCP will be no larger. In case (iii), it is the same as the original. In the case where there is a class for which the abandonment rate is zero, the relation between abandonment and queue lengths cannot be used. However, this can be treated similarly to the following paragraph.

Consider next the case where $\ell$ is an increasing function of $\Phi_i$, for all $i$, with $\theta_i \leq \mu_i$ for all $i$. Here, when a customer is moved into service, its service time is coupled to the abandonment time of its twin: Pick an exponential random variable with rate $\mu_i$, and a Bernoulli random variable that is 1 with probability $\theta_i/\mu_i$. The service time is the exponential random variable. If the Bernoulli random variable is 1, then the abandonment occurs simultaneously; otherwise the original customer does not abandon at that time and picks a new exponential random variable with rate $\theta_i$. Again, if the original SCP needs the server, the customer is moved out. It can be seen that the cost of the modified SCP will be no larger than the original one.

5.2. *Additional topics.* The following is a list of research problems that are suggested by the present study.

1. *Nonlinear waiting costs*: Nonlinear waiting costs are natural for quantifying human costs of waiting [37, 39]. We believe that it is possible to reduce such costs to nonlinear costs of queue lengths, and are planning to include this in future work.
2. *Alternative cost structures*: Discounted costs are mathematically convenient. Long-run average costs provide an alternative which is no less, perhaps more, natural for call center applications. Their analysis, however, would be mathematically more taxing.
3. *Performance analysis in the QED regime*: In the present study, we are not analyzing the performance of our queueing system under the proposed SCPs. In particular, one would like to confirm that the (discounted) probability of delay, for each class, is nontrivial, as expected in the QED regime. Such analysis might require numerical supplements, as in [21]. This could also shed further light on qualitative features of our asymptotically optimal SCPs.
4. *More general models*: The model in Figure 1 is a beginning. Ultimately, one would like to generalize it to the model surveyed in [38], which has heterogeneous pools of servers with overlapping service skills. (See [14] for interesting simulations of such models.) In conventional heavy traffic (efficiency driven), a simple generalized $C\mu$ control was proved asymptotically optimal [30]. Here, only the problem of assigning servers who become idle is relevant, since customers essentially never encounter an idle server upon arrival. This same simplifying feature applies for our model, under work conservation. But with heterogeneous pools of servers, and with a nontrivial fraction of arrivals encountering idle servers (as expected in the QED regime), both the assignment of servers to customers and the routing of arriving customers to idle servers become significant. In a call center context, the problem of online matching customers and servers is called skills-based routing; it is widely acknowledged as the most important and difficult operational problem next to staffing, to which we now turn.

5. *Staffing insights*: The staffing problem is to determine the least (optimal) number of servers $n$ that is required to conform to given performance standards. In the QED regime, $n \approx R + \beta\sqrt{R}$, where $R$ is the offered load and $\beta$ is a scalar. The problem can thus be decomposed, as in [7], into two steps: first, given a QED operation, determine the least (optimal) scalar $\beta$; then, establish that operating in the QED regime is indeed desirable (optimal). The staffing problem becomes more interesting and far more difficult in a skills-based routing environment. ([8] is the single paper on the subject that we are aware of.)

## APPENDIX

PROOF OF PROPOSITION 2. Note that $(x, u) \mapsto b(x, u)$ is continuous and $x \mapsto b(x, u)$ is Lipschitz uniformly in $u$. Consider $b_m$, a function that agrees with $b$ on the ball $B(0, m)$, uniformly Lipschitz and bounded. Then strong existence and uniqueness for

$$X_m(t) = x + rW(t) + \int_0^t b_m(X_m(s), u(s))\,ds, \qquad 0 \le t < \infty,$$

holds by Theorem I.1.1 of [6]. Since $\|X_m(t)\| \le \|x\| + c\|W(t)\| + c\int_0^t \|X_m(s)\|\,ds$, one has $\|X_m(t)\| \le (\|x\| + c\|W\|_t^*)(1 + e^{ct})$ by Gronwall's lemma. Thus letting $\tau_m = \inf\{t : \|X_m(t)\| \ge m\}$, one has $\tau_m \to \infty$ a.s. Therefore $X(t) = \lim_m X_m(t)$ for all $t$ defines a process that solves the equation (a strong solution). If $X$ and $\bar{X}$ are both strong solutions, then, for every $m$, they both agree with $X_m$ on $[0, \tau_m]$. Therefore they agree on $[0, \infty)$ a.s. $\square$

PROOF OF LEMMA 1. Let $n \in \mathbb{N}$ and $i \in K$ be fixed, and consider for each $s$ the $\sigma$-fields

$$\bar{F}_s = \sigma\{\mathbf{1}_{\{\mathring{T}_i^n(u) \le s\}}, R_i^n(\alpha) : u \in \mathbb{R}_+, \alpha \le s\},$$

$$\bar{G}_s = \sigma\{R_i^n(\beta + \gamma) - R_i^n(\beta) : \beta > s, \gamma > 0\}.$$

We simplify notation by writing $T_u = \mathring{T}_i^n(u)$ and $R(u) = R_i^n(u)$.

For each $t$ and $s$, one has $\{T_t \le s\} \in \bar{F}_s$, and therefore, for each $t$, $T_t$ is a stopping time on the filtration $(\bar{F}_s)$. We next show that $M_t := R(t) - \theta_i t$ is a martingale on the filtration $(\bar{F}_s)$; hence the lemma follows from the optional stopping theorem. Indeed, it is clear that $M_s$ is measurable on $\bar{F}_s$ for each $s$. Moreover, $M_r - M_s$ is measurable on $\bar{G}_s$ for each $s$ and $r \geq s$. It remains to show that $\bar{F}_s$ is independent of $\bar{G}_s$ for each $s$. Fix $s$. Fix $\delta$, and $u \geq 0$, $0 < \alpha < s < s + K\delta = \beta$, $\gamma > 0$. Let

$$H_m = \{T_{m\delta} \leq s < T_{(m+1)\delta}\}.$$

Note that $P(\bigcup_m H_m) = 1$. Let

$$\bar{H}_{K,r} = \{T_r > s;\ \dot{T} \leq K \text{ on } [0,r]\}.$$

Let also

$$\hat{H}_{m,K} = \{\dot{T} \leq K,\ \text{on } [m\delta, (m+1)\delta]\}.$$

For measurable bounded $f, g$ (denote by $c$ a bound on $fg$),

$$C_{f,g} := E[f(\mathbf{1}_{\{T_u \leq s\}}, R(\alpha))g(R(\beta+\gamma) - R(\beta))]$$

$$= \sum_{m=0}^{[r/\delta]} E[\mathbf{1}_{H_m \cap \hat{H}_{m,K}} f(\mathbf{1}_{\{T_u \leq s\}}, R(\alpha))g(R(\beta+\gamma) - R(\beta))] + e_1,$$

where $|e_1| \leq cP(\bar{H}_{K,r}^c)$. Under the event $H_m \cap \hat{H}_{m,K}$, $T_{m\delta} \leq s \leq T_{(m+1)\delta} \leq s + K\delta = \beta$. Denote $\Delta_m = \beta - T_{(m+1)\delta}$, and note that $0 \leq \Delta_m \leq K\delta$ under the same event. Then

$$C_{f,g} = \sum_{m=0}^{[r/\delta]} E\{\mathbf{1}_{H_m \cap \hat{H}_{m,K}} f(\mathbf{1}_{\{T_u \leq s\}}, R(\alpha))g(R(T_{(m+1)\delta} + \Delta_m + \gamma)$$
$$- R(T_{(m+1)\delta} + \Delta_m))\} + e_1.$$

Let $\tilde{H}_{K,\delta}$ denote the event that there are no jumps of the process $R$ within $[s, s+K\delta] \cup [s+\gamma, s+\gamma+K\delta]$:

$$C_{f,g} = \sum_{m=0}^{[r/\delta]} E\{\mathbf{1}_{H_m \cap \hat{H}_{m,K} \cap \tilde{H}_{K,\delta}} f(\mathbf{1}_{\{T_u \leq s\}}, R(\alpha))g(R(T_{(m+1)\delta} + \gamma) - R(T_{(m+1)\delta}))\}$$
$$+ e_1 + e_2$$
$$= \sum_{m=0}^{[r/\delta]} E\{\mathbf{1}_{H_m \cap \hat{H}_{m,K}} f(\mathbf{1}_{\{T_u \leq s\}}, R(\alpha))g(R(T_{(m+1)\delta} + \gamma) - R(T_{(m+1)\delta}))\}$$
$$+ e_1 + e_2 + e_3,$$

where $|e_2|, |e_3| \leq cP(\tilde{H}_{K,\delta}^c)$. Recall that by Definition 2(i), $\mathcal{F}^n(t)$ and $\mathcal{G}^n(t)$ are independent. Since under $H_m$, $\alpha \leq s \leq T_{(m+1)\delta}$, it follows that $\mathbf{1}_{H_m \cap \hat{H}_{m,K}} f(\mathbf{1}_{\{T_u \leq s\}}, R(\alpha)) \in \mathcal{F}^n((m+1)\delta)$. Also, $R(T_{(m+1)\delta} + \gamma) - R(T_{(m+1)\delta}) \in \mathcal{G}^n((m+1)\delta)$, and, using Definition 2(ii), it has the same law as $R(\gamma)$. Hence

$$C_{f,g} = Eg(R(\gamma)) \sum_{m=0}^{[r/\delta]} E\{\mathbf{1}_{H_m \cap \hat{H}_{m,K}} f(\mathbf{1}_{\{T_u \le s\}}, R(\alpha))\} + e_1 + e_2 + e_3$$

$$= Eg(R(\gamma)) \sum_{m=0}^{\infty} E\{\mathbf{1}_{H_m} f(\mathbf{1}_{\{T_u \le s\}}, R(\alpha))\} + e_1 + e_2 + e_3 + e_4$$

$$= Eg(R(\gamma)) E\{f(\mathbf{1}_{\{T_u \le s\}}, R(\alpha))\} + e_1 + e_2 + e_3 + e_4,$$

where $|e_4| \le cP(\bar{H}_{K,r}^c)$. It follows from Lemma 2 that $E\|\Phi^n\|_r^* < \infty$. Since $\dot{T} = \Phi_i^n$, one has that

$$\lim_{r \to \infty} \liminf_{K \to \infty} P(\bar{H}_{K,r}) = 1.$$

Note also that $P(\tilde{H}_{K,\delta}^c) \le c_1 K \delta$ for some constant $c_1$. Taking $\delta \to 0$ and $K \to \infty$ such that $K\delta \to 0$, and then taking $r \to \infty$, we conclude that

$$E[f(\mathbf{1}_{\{T_u \le s\}}, R(\alpha))g(R(s+\gamma) - R(s))]$$
$$= E[f(\mathbf{1}_{\{T_u \le s\}}, R(\alpha))]E[g(R(\gamma))]$$
$$= E[f(\mathbf{1}_{\{T_u \le s\}}, R(\alpha))]E[g(R(s+\gamma) - R(s))].$$

Since $\alpha, s$ and $\gamma$ are arbitrary (subject to $0 < \alpha < s < s+\gamma$), and so are $f, g$, it follows that $\bar{F}_s$ and $\bar{G}_s$ are independent for any $s$. The result follows. $\square$

PROOF OF PROPOSITION 1 (SKETCH).

Existence and uniqueness for the system (7) and (11) are easily obtained by induction on the jump times of the processes $A_i^n$, $R_i^n$ and $S_i^n$. By the assumptions on the function $F$, the constraints (6) are met. We next need to show that Definition 2 holds. For part (i) of the definition it suffices to show that, for any bounded measurable $g$,

$$E[g(A_i^n(\tau_i^n(t) + u) - A_i^n(\tau_i^n(t)), \ S_i^n(T_i^n(t) + u) - S_i^n(T_i^n(t)),$$
(84) $$R_i^n(\mathring{T}_i^n(t) + u) - R_i^n(\mathring{T}_i^n(t)); i \in K)|\mathcal{F}_t^n]$$
$$= Eg(A_i^n(u), S_i^n(u), R_i^n(u); i \in K)$$

where $u > 0$, and for part (ii) is suffices to show that

$$E[g(S_i^n(T_i^n(t) + u_j) - S_i^n(T_i^n(t)),$$
(85) $$R_i^n(\mathring{T}_i^n(t) + u_j) - R_i^n(\mathring{T}_i^n(t)); i \in K, \ j \ge 1)|\mathcal{F}_t^n]$$
$$= Eg(S_i^n(u_j), R_i^n(u_j); i \in K, j \ge 1),$$

where $0 < u_1 < u_2 < \cdots$. In what follows we suppress $n$ from the notation, fix $i$ and $t$ and, denoting $T_t^i = T_i^n(t)$, show that, for $u > v > 0$,

(86) $$E[g(S_i(T_t^i + u) - S_i(T_t^i + v))|\mathcal{F}_t] = Eg(S_i(u - v)),$$

(87) $$E[g(A_i(\tau_i(t) + u) - A_i(\tau_i(t) + v))|\mathcal{F}_t] = Eg(A_i(u - v)).$$

Since the notation is quite complicated, we do not give the full details on proving (84), (85), but only comment that the argument is similar to the one we use in proving (86) and (87).

To show (86), for $\delta > 0$, let $H_m = \{T_t^i \in [m\delta, (m+1)\delta)\}$ and

$\tilde{H}_m = \{S_i \text{ has no jumps on } [m\delta + v, (m+1)\delta + v] \cup [m\delta + u, (m+1)\delta + u]\}$.

Then

$$E[g(S_i(T_t^i + u) - S_i(T_t^i + v))|\mathcal{F}_t]$$
$$= \sum_{m=0}^{\infty} E[\mathbf{1}_{H_m} g(S_i(T_t^i + u) - S_i(T_t^i + v))|\mathcal{F}_t]$$
$$= \sum_{m=0}^{\infty} E[\mathbf{1}_{H_m} g(S_i((m+1)\delta + u) - S_i((m+1)\delta + v))|\mathcal{F}_t] + e_1$$
$$= \sum_{m=0}^{\infty} \mathbf{1}_{H_m} E[g(S_i((m+1)\delta + u) - S_i((m+1)\delta + v))|\mathcal{F}_t] + e_1,$$

where

(88) $$|e_1| \leq c \sum_{m=0}^{\infty} P(H_m \cap \tilde{H}_m^c | \mathcal{F}_t).$$

Note that on the event $T_t^i \leq \eta$, the quantities $(X(s), \Psi(s), \Phi(s); s \leq t)$ only depend on $A$, $R$, $S_j$, $j \neq i$, and $S_i(t')$, $t' \leq \eta$. Since $S_i$ is Poisson and independent of the processes $A$, $R$ and $S_j$, $j \neq i$, using the definition of $\mathcal{F}_t$ and $H_m$ we obtain that

$$E[g(S_i(T_t^i + u) - S_i(T_t^i))|\mathcal{F}_t]$$
$$= \sum_{m=0}^{\infty} \mathbf{1}_{H_m} Eg(S_i((m+1)\delta + u) - S_i((m+1)\delta)) + e_1$$
$$= Eg(S_i(u)) + e_1.$$

By (88), and since $\tilde{H}_m$ depends only on $S_i(s)$; $s \geq m\delta + v$,

$$|e_1| \leq c \sum_{m=0}^{\infty} \mathbf{1}_{H_m} P(\tilde{H}_m^c) \leq c\delta,$$

where $c$ does not depend on $\delta \in (0,1)$. As a result, (86) holds. An equivalent of (86) for the processes $R_i$ is proved analogously. Equality (87) is proved analogously, where one conditions on $\mathcal{F}_t \vee \sigma\{\tau_i(t)\}$ and uses the fact that $A_i$ is a renewal process. $\square$

PROOF OF PROPOSITION 3. Throughout, fix a compact subset $A$ of $\mathbb{X}$, and let $c$ denote a positive constant that depends only on $A$, and whose value may change from location to location. Recall from the proof of Theorem 3 that, for each $x \in \mathbb{R}^k$, $h(x)$ satisfies $\varphi(x, h(x)) = \inf_{v \in \mathbb{S}^k} \varphi(x, v)$, where

$$\varphi(x, u) = b(x, u) \cdot Df(x) + L(x, u).$$

In the special case we analyze here, $L(x, u) = \sum_i g_i((\mathbb{1} \cdot x)^+ u_i)$, hence [cf. (26)]

$$\varphi(x, u) = (\ell + (\mu - \theta)(\mathbb{1} \cdot x)^+ u - \mu x) \cdot Df(x) + \sum_i g_i((\mathbb{1} \cdot x)^+ u_i)$$

$$=: \bar{a}(x) + \bar{b}(x) \cdot u + \sum_i g_i(\bar{x} u_i),$$

where $\bar{x} = \mathbb{1} \cdot x > 0$. For any $x \in \mathbb{X}$, the map $u \mapsto \varphi(x, u)$ is strictly convex; hence the infimum over $\mathbb{S}^k$ is uniquely attained.

Fixing $x \in A$, and letting $m_i(u_i) = \bar{b}_i(x) u_i + g_i(\bar{x} u_i)$, $\varphi(x, u)$ is given as $\bar{a}(x) + \sum_i m_i(u_i)$. Use Taylor's formula for each $m_i$ based at $u_i$,

$$\varphi(x, v) = \varphi(x, u) + \sum_i m_i'(u_i)(v_i - u_i) + (1/2) m_i''(\xi_i)(v_i - u_i)^2.$$

We claim that $\sum_i m_i'(u_i)(v_i - u_i) \geq 0$ for $v \in \mathbb{S}^k$. For if this is false, let $v \in \mathbb{S}^k$ be such that $\sum_i m_i'(u_i)(v_i - u_i) = -c < 0$. Then for $v^\varepsilon := u + \varepsilon(v - u)$, $\sum_i m_i'(u_i)(v_i^\varepsilon - u_i) = -c\varepsilon$. Moreover, by assumption on the functions $g_i$, there is a constant $c$ such that $g_i''(\bar{x} v_i) \leq c$; hence $m_i''(v_i) \leq \bar{x}^2 c$, $v \in \mathbb{S}^k$. Therefore $\sum_i m_i''(\xi_i)(v_i^\varepsilon - u_i)^2 \leq c\varepsilon^2$, implying that $\varphi(x, v^\varepsilon) < \varphi(x, u)$ for $\varepsilon > 0$ small, contradicting the definition of $u$.

Using the above, and that $m_i''(\xi_i) = \bar{x}^2 g_i''(\bar{x} \xi_i) \geq c_0 \bar{x}^2 \geq c > 0$ on $A$, we obtain

$$\begin{aligned}(89) \quad \varphi(x, v) - \varphi(x, u) &\geq (1/2) \sum_i m_i''(\xi_i)(v_i - u_i)^2 \\ &\geq c \|v - u\|^2, \qquad x \in A, \ u = h(x), \ v \in \mathbb{S}^k.\end{aligned}$$

Let $x, y \in A$ and let $u = h(x)$ and $v = h(y)$. Since $f$ is of class $C^2$ (cf. Theorem 3),

$$\|Df(x) - Df(y)\| \leq c\|x - y\|.$$

By the proof of Theorem 3,

$$|H(x, p) - H(y, q)| \leq c(\|p - q\| + \|x - y\|^\varrho),$$

for $p, q$ in a compact set. It follows that

$$|\varphi(x,u) - \varphi(y,v)| = |H(x, DV(x)) - H(y, DV(y))|$$
$$\leq c\|x-y\|^\varrho.$$

Since by Assumption 2(iii) on $L$, $x \mapsto \varphi(x,v)$ is Hölder of exponent $\varrho$,

$$\varphi(x,v) - \varphi(x,u) \leq c\|x-y\|^\varrho.$$

Combining the last display with (89), $\|u-v\|^2 = \|h(x) - h(y)\|^2 \leq c\|x-y\|^\varrho$, and the result follows. $\square$

**Acknowledgments.** The authors are grateful to Paul Dupuis for valuable discussions, to Michael Harrison and Assaf Zeevi for access to an early version of [21], and to an associate editor and two referees for their thoughtful comments.

R. ATAR
DEPARTMENT OF ELECTRICAL
  ENGINEERING
TECHNION–ISRAEL INSTITUTE OF
  TECHNOLOGY
HAIFA 32000
ISRAEL
E-MAIL: atar@ee.technion.ac.il

A. MANDELBAUM
DEPARTMENT OF INDUSTRIAL
  ENGINEERING AND MANAGEMENT
TECHNION–ISRAEL INSTITUTE OF
  TECHNOLOGY
HAIFA 32000
ISRAEL

M. I. REIMAN
BELL LABS, LUCENT TECHNOLOGIES
600 MOUNTAIN AVENUE
MURRAY HILL, NEW JERSEY 07974
USA